\documentclass[12pt]{article}
\usepackage{amsmath,epsfig, amsthm,amsfonts,amsbsy,amssymb,latexsym,amsxtra}

\makeatletter \@addtoreset{equation}{section}

\newcommand{\beq}[1]{\begin{equation} \label{#1}}
\newcommand{\eeq}{\end{equation}}
\newcommand{\bed}{\begin{displaymath}}
\newcommand{\eed}{\end{displaymath}}
\newcommand{\ben}{\begin{eqnarray*}}
\newcommand{\een}{\end{eqnarray*}}

\def\F{{\cal F}}
\def\op{{\cal L}}
\def\cd{(\cdot)}
\def\M{{\cal M}}
\def\rr{{\mathbb R}}

\newcommand{\la}{\lambda}

\newcommand{\e}{\varepsilon}

\newcommand{\al}{\alpha}

\newcommand{\sg}{\sigma}
\newcommand{\tha}{\theta}
\newcommand{\pr}{{\mathbf P}}
\newcommand{\ex}{{\mathbf E}}
\newcommand{\wdt}{\widetilde}

\newcommand{\xalz}{\ensuremath{x_0,\alpha_0}}
\newcommand{\xz}{\ensuremath{x_0}}

\newtheorem{thm}{Theorem}[section]
\newtheorem{prop}[thm]{Proposition}
\newtheorem{lem}[thm]{Lemma}

\theoremstyle{definition}
\newtheorem{rem}[thm]{Remark}
\newtheorem{defn}[thm]{Definition}
\newtheorem{exm}[thm]{Example}

\newcommand{\thmref}[1]{Theorem~{\rm \ref{#1}}}
\newcommand{\lemref}[1]{Lemma~{\rm \ref{#1}}}

\newcommand{\propref}[1]{Proposition~{\rm \ref{#1}}}

\newcommand{\set}[1]{\left\{#1\right\}}
\newcommand{\abs}[1]{\left\vert #1\right\vert}

\newcommand{\disp}{\displaystyle}

\def\({\left(}
\def\){\right)}

\def\one{{\hbox{1{\kern -0.35em}1}}}

\newcommand{\bea}{\bed\begin{array}{rl}}
\newcommand{\eea}{\end{array}\eed}
\newcommand{\ad}{&\!\!\!\disp}
\newcommand{\aad}{&\disp}
\newcommand{\barray}{\begin{array}{ll}}
\newcommand{\earray}{\end{array}}
\newcommand{\lbar}{\overline}

\title{On Optimal Harvesting Problems in Random Environments}
\author{Qingshuo Song\thanks{Department of Mathematics, City
    University of Hong Kong, 83 Tat Chee Avenue, Kowloon Tong, Hong
    Kong, {\tt   song.qingshuo@cityu.edu.hk}. The research of this author was supported in part by the Research Grants Council of Hong Kong No. CityU 100310. }
    \and Richard H. Stockbridge\thanks{Department of Mathematical Sciences, University of
    Wisconsin-Milwaukee, Milwaukee, WI 53201, {\tt  stockbri@uwm.edu}.
    The research of this author was supported in part by the U.S. National Security Agency under Grant Agreement Number H98230-09-1-0002.}
     \and Chao   Zhu\thanks{Department of Mathematical Sciences, University of
    Wisconsin-Milwaukee, Milwaukee, WI 53201, {\tt zhu@uwm.edu}.} }
\begin{document}

\maketitle

\begin{abstract}
This paper investigates the optimal harvesting strategy for a single species living in random environments whose growth is given by a regime-switching diffusion.  Harvesting acts as a (stochastic) control on the size of the population.  The objective is to find a harvesting strategy which maximizes the expected total discounted income from harvesting {\em up to the time of extinction}\/ of the species; the income rate is allowed to be state- and environment-dependent.  This is a singular stochastic control problem with both the extinction time and the optimal harvesting policy depending on the initial condition.  One aspect of receiving payments up to the random time of extinction is that small changes in the initial population size may significantly alter the extinction time when using the same harvesting policy.  Consequently, one no longer obtains continuity of the value function using standard arguments for either regular or singular control problems having a fixed time horizon.  This paper introduces a new sufficient condition under which the continuity of the value function for the regime-switching model is established.  Further, it is shown that the value function is a viscosity solution of a coupled system of quasi-variational inequalities.  The paper also establishes a verification theorem and, based on this theorem, an $\varepsilon$-optimal harvesting strategy is constructed under certain conditions on the model.  Two examples are analyzed in detail.



\end{abstract}

\setlength{\baselineskip}{0.21in}
\section{Introduction}\label{sect-intro}

One of the most important yet difficult problems in modern natural resources management is the establishment of ecologically, environmentally and economically reasonable wildlife management and  harvesting policies.  There are many occurrences where myopic unconstrained harvesting has led to  local and/or global extinctions. Lande et al. \cite{Lande} documents many such examples.  Real ecological communities are random by nature.  As a result of developments in stochastic analysis and stochastic control techniques, there has been a resurgent interest in determining the optimal harvesting strategies in the presence of stochastic fluctuations (see, e.g.,
\cite{A-Shepp, Brauman-02,Lungu-O,Ryan-Hanson}). Unfortunately, most of the current research on harvesting problems,  including the aforementioned references, are primarily  focused on a single species in a static environment.  The paper by Lungu and {\O}ksendal \cite{L-Oksendal} makes a first step in the analysis of the harvesting problem for interacting populations but does not consider changes in the environment.

As noted in \cite{Cohen-90,Jeffries76,Slatkin-78}, the variations in the
external environment (for example, weather or anthropogenic) can have important
 effects on the dynamics of the populations of the ecosystem.  In addition to the
  random fluctuations of the populations (usually  modeled by white noise; see, e.g., %Arnold et al.
  \cite{Arnold-H-S-79}),
   certain biological parameters such as the growth rates and the carrying capacities
    often demonstrate abrupt changes due to environmental noise.  Moreover,
     the qualitative changes of those parameters form an essential part of the
     dynamics of the ecosystem. For example, Medina-Reyna \cite{Medina-R} demonstrates
     that the mean growth rates of white shrimp
     (Litopenaeus vannamei) in the Mar Muerto Lagoon, Southern Mexico
      are significantly different in various salinity levels.  Similar observations were made in
   %Otuma and Osakwe
   \cite{Otuma-O}
   for reproduction performance of crossbred goats in a derived Guinea savanna zone. For another example, the carrying capacities often vary according to the changes in nutrition, water supply, living spaces, and/or food resources (see
   %Sayre
   \cite{Sayre} for many such examples).  In the mathematical community, people are paying more attention to the modeling and analysis of population dynamics subject to both white and colored noises; see, e.g., \cite{Du-Sam-06,Luo-Mao07,ZY-09a} and references therein.

Naturally, one expects that the optimal harvesting strategies may vary according to the changes of the environment.  Despite the increasing interests in the mathematical modeling and analysis of population dynamics, to our best knowledge, there are relatively few results in the literature that address harvesting strategies in random environments.  This paper addresses this hole in the literature by examining optimal harvesting problems of a species in random environments.

Suppose there is  a single  species whose growth is subject to the usual fluctuations as well as the abrupt changes of the random environments.  Harvesting strategies are introduced to derive financial benefit as well as to control the growth of the population.  The goal is to find a harvesting strategy which maximizes the expected total discounted income from harvesting, up to the time when the  population falls to a given threshold (e.g., extinction).  Harvesting may occur instantaneously so results in a singular stochastic control problem in the sense that the optimal harvesting strategy may not be absolutely continuous with respect to the Lebesgue measure of time. In other words, in contrast to the regular stochastic control problems, in which the displacement of the state due to control is differentiable in time, the harvesting problem considered in this work allows the displacement to be discontinuous.  This paper establishes a verification theorem and, based on the theorem, explicitly constructs an $\e$-optimal harvesting strategy.  Both the extinction time and harvesting policy may depend on the initial conditions.  As a result, continuity of the value function can not be obtained using the standard arguments for regular or singular stochastic control problems in a fixed time horizon.  This paper provides a sufficient condition under which the continuity of the value function is guaranteed. It is further shown that the value function is a viscosity solution of a coupled system of  quasi-variational inequalities \eqref{quasi-vari}.

% The novelty  of this work is the consideration of abrupt changes
% in the random environments, which significantly affects the dynamics
% of the population size of the species and hence the optimal
% harvesting strategy (see the discussions in Example \ref{example-2regimes}
% for details).  Another important feature  is the
% state-and-regime-dependent price function $f$ (defined in Section
% \ref{sect-formulation}), reflecting the fact that in many practical
% situations, the marginal revenues from harvesting depend on both the
% current state of the system and the current regime of the environment.
% While these considerations allow a better approximation of the real world dynamics,
% they also introduce greater difficulty for the analysis.

The novelty of this work arises in two distinct ways.  The modeling of random environments through the use of a continuous-time finite-state Markov chain introduces coupling of the value function for each environment in the quasi-variational inequalities for the verification theorem (Theorem \ref{thm-verification}).  An $\varepsilon$-optimal harvesting policy (Theorem \ref{cor-1d-optimal}) is determined under certain conditions that involves quickly harvesting very small amounts until the species becomes extinct in a very small time interval. We should also remark that the proof of Theorem \ref{cor-1d-optimal}  is
    very technical and non-trivial. In addition
to the subtle analysis in dealing with the controlled process $\hat X$,
the presence of environmental switching adds  much difficulties in the proof.  The introduction of different regimes necessarily implies that the optimal harvesting strategy will depend on the current environment (Example \ref{example-2regimes}).  The determination of an environment-dependent optimal policy is non-trivial as one must overcome some significant technical challenges.

This paper's second contribution comes from identifying a new sufficiency condition for the continuity of the value function as a function of the initial state.  The fact that payment is received only until the random time of extinction and this time strongly depends on the harvesting policy adopted means that a small decrease in initial population size may result in a significant decrease of the extinction time.  Continuity is therefore not a direct extension of standard results for a fixed time horizon.  Theorem \ref{thm-value-cont} establishes a sufficient condition under which the value function can be proven to be continuous for this criterion involving the hitting time of the population at $0$ (or any quasi-extinction level from which the population will not rebound).  Once the value function is shown to be continuous, it is then proven to be a viscosity solution of the quasi-variational inequalities (Theorem \ref{thm-viscosity}); even this analysis is technically challenging due to the existence of multiple environments.

Note that this work is expressed entirely in terms of harvesting of a single species in random environments, but as in Miller and Voltaire \cite{Miller-V}, the harvesting problem is a paradigm that has many additional economic applications.

Besides the optimal harvesting problems considered in this paper and  \cite{A-Shepp,L-Oksendal,Lungu-O}, singular stochastic control has found applications in many other areas.  For example, singular stochastic control problems naturally arise in monotone follower problems \cite{Karatzas-S}, optimal dividend distribution schemes \cite{Choulli,Paulsen-03}, portfolio selection management with transaction cost \cite{MSXZ08a,OS02}, diffusion control of many-server queues \cite{Weerasinghe}, and heavy traffic modeling and control problems \cite{Wein90}.
We refer the reader to \cite{FlemingS,Kushner-D,YongZ} for more such examples.  See also \cite{Hauss-S-I,Hauss-S} for  a general singular stochastic control problem for a multidimensional It\^o diffusion on a fixed time horizon,
in which
the existence of the optimal control and the characterization of the value function as the unique viscosity solution of a Hamilton-Jacobi-Bellman equation are established.
Most, if not all, of the existing literature on singular stochastic controls consider with It\^o (jump) diffusions.

Regular control and optimal stopping problems for regime switching diffusions have become more popular recently (see, for example, \cite{Guo-Z,Guo-Z05,Taksar-Z,Zhou-Yin} and references therein).  Less is known for singular control of regime switching diffusions.  Moreover, a common assumption is that the marginal yield from exerting the singular control is constant.  Two exceptions are in \cite{Alvarez,Lungu-O}, where the marginal yields depend  on state and time, respectively.  The assumption of constant marginal yield seems rather restrictive since in the real world, the unit price usually depends on the current state of the system.  This paper considers state-and-regime-dependent marginal yields from harvesting; that is, the unit price depends on both the current state of the population size of the species and the regime of the environment.  This additional feature of the model is not merely an extension of the traditional models, but in fact, introduces many interesting mathematical problems for the analysis.  More specifically,
there may not exist {\em admissible} optimal harvesting strategy  under this setting; see \thmref{cor-1d-optimal}, Remark \ref{rem-optimal-control}, and Example \ref{example-non-constant-price} for more details. Nevertheless, using detailed and careful analysis of the sample path properties of the controlled process, we constructs an explicit admissible $\e$-optimal harvesting strategy.

The rest of the paper is organized as follows. A precise formulation of the problem is presented in Section \ref{sect-formulation}.  Then a verification theorem is proven in Section \ref{sect-verification} and is used to explicitly construct an $\e$-optimal harvesting strategy under additional conditions.  Two examples are given in Section \ref{sect-example} to illustrate these results.  Section \ref{sect-properties-of-V} derives the continuity of the value function $V$ and characterizes it as a viscosity solution of a coupled system of quasi-variational inequalities \eqref{quasi-vari}.  Section \ref{sect-conclusion} contains concluding remarks.

A few words about notation is needed. A function from $[0,\infty)$ to some Polish space $E$ is {\em c\`adl\`ag} if it is right continuous and has left limits in $E$. When $E=\rr$ and $\xi $ is c\`adl\`ag, then $\Delta \xi(t)= \xi(t)-\xi(t-)$ for $t> 0$ and the convention $\Delta \xi(0)= \xi(0)$ is used. As usual, $\sup \emptyset=-\infty$ and $\inf\emptyset= + \infty$.  For any $a,b \in \rr$, $a^+ =\max\set{a,0}$ and $a\wedge b= \min \set{a,b}$. If $B $ is a set, $I_B$ denotes the indicator function of $B$.

\subsection{Formulation}\label{sect-formulation}

Suppose a certain species, whose population size at time $t$ is denoted by $X(t)$, lives in random environments.  As alluded in Section \ref{sect-intro}, in addition to the  random fluctuations of the population, we also assume that the growth of the species is subject to abrupt changes of the environment.
%   (usually  modeled by white noise, see, for example, %Arnold et al.
%   \cite{Arnold-H-S-79}),   certain biological parameters such as the growth rates and the carrying capacities
%    often demonstrate abrupt changes due to environmental noise.
%    Moreover, the qualitative changes of those parameters form an essential part of the dynamics of the ecosystem. For example,
%   %the growth rates of some species in the rainy season will be much different from those in the dry season.
%    Medina-Reyna \cite{Medina-R} demonstrates that the mean growth rates of white shrimp (Litopenaeus vannamei)
%    in the Mar Muerto Lagoon, Southern Mexico are significantly different in different salinity levels.
%     Similar observations were made in %Otuma and Osakwe
%      \cite{Otuma-O} for reproduction performance of
%     crossbred goats in a derived Guinea savanna zone. For another example, the carrying capacities often vary according
%     to the changes in nutrition, water supply, living spaces, and/or food resources, see
%    %Sayre
%    \cite{Sayre} for many such examples.  We refer to %Cohen et al.
%     \cite{Cohen-90,Jeffries76,Slatkin-78} %, Jeffries \cite{Jeffries76},   Slatkin \cite{Slatkin-78},
%      and references therein for more details and examples.
For simplicity, we assume that the switching among different environments is memoryless and the waiting time for the next switch is exponentially distributed.  In fact, this phenomenon is also frequently observed in the nature; see the aforementioned references. Thus we can model the random environments and other random factors in the ecological system by a continuous-time Markov chain $\set{\al(t), t\ge 0}$ with a finite state space $\M=\set{1, \dots,m}$.  Let the continuous time Markov chain $\al\cd$ be generated by $Q=(q_{ij})$, that is,
  \begin{equation}\label{Q-gen}\pr\set{\al(t+ \Delta t)=j|
\al(t)=i,\al(s),s\le t}=\begin{cases}q_{ij}
\Delta t + o(\Delta t),\ &\hbox{ if }\ j\not= i\\
1+ q_{ii}\Delta t + o(\Delta t), \ &\hbox{ if }\ j=i,
\end{cases}   \end{equation}  where $q_{ij}\ge 0$ for $i,j=1,\dots,m$
with $j\not= i$ and $ q_{ii}=-\sum_{j\not= i}q_{ij}<0$ for each $i=1,\dots,m$.

 In light of the above discussion, in an effort  to capture the salient  feature that continuous dynamics and discrete events coexist in the ecosystem, we model the evolution of $X (t)$ {\em in the absence of harvesting}   by the stochastic differential equation
 \begin{equation}\label{dyna-compo} dX (t) =b (X(t),\al(t))dt + \sigma(X(t),\al(t))dw(t), \ \ X(0)=x,\al(0)=\al, \end{equation}
 where  $w(\cdot) $ is a $1$-dimensional standard Brownian motion which provides the random fluctuations in the population's size,
 %            $\al\cd \in \M=\set{1,\dots,m}$ is a continuous time Markov chain modeling the changes of the random environments,
  and   $b$ and $\sigma$ are real-valued functions.
   Further, we assume that the Brownian motion $w\cd$ and the Markov chain $\al\cd$ are independent, a standard assumption in the literature.

 Assume throughout the paper that
 $b$ and $\sigma$ satisfy the usual Lipschitz condition and the linear growth condition. That is, there exists some $\kappa_0 >0$ such that for any $x,y \in \rr $ and each $\al \in \M$, we have
  \begin{equation}\label{ito-condition}\barray \ad \abs{b(x,\al)-b(y,\al)} + \abs{\sigma(x,\al)-\sigma(y,\al)} \le \kappa_0 \abs{x-y},  \\ \ad \abs{b(x,\al)} + \abs{\sigma(x,\al)} \le \kappa_0(1+ \abs{x}).
  \earray \end{equation}
Consequently, the solution $X^{x,\al}\cd$ of (\ref{dyna-compo}) exists and is   unique in the strong sense (see \cite{YZ-10}  for details).
   Moreover, the solution  $X^{x,\al}\cd$ will not explode in finite time with probability 1 or it is {\em regular} in the sense of Khasminskii \cite{K}.
 We refer the reader to \cite{YZ-10} for related results on the regularity of regime switching diffusions.

 If the species is subject to harvesting, and if $Z(t)$ denotes the total amount harvested from the   species up to time $t$, then $\hat X\cd $,  the population size of the {\em harvested population}, satisfies
  \begin{equation}\label{harvest} d\hat X(t) =b(\hat X(t),\al(t)) dt + \sigma(\hat X(t),\al(t))dw(t)-dZ(t),\end{equation}
 with initial conditions
 \begin{equation}\label{initial} \hat X(0-) = x \in \rr_+,\ \  \al(0)= \al \in \M.  \end{equation}
Note that $\hat X(0)$ may not equal to $\hat X(0-)$ due to an instantaneous harvest $Z(0)$ at time $0$.
Throughout the paper we use the convention  that $Z(0-)=0$. The jump size of $Z$ at time $t\ge 0$ is denoted by $\Delta Z(t):= Z(t)-Z(t-) $, and $Z^c(t) := Z(t)-\sum_{0\le s \le t} \Delta Z(s)$ denotes the continuous part of $Z$.
Also note that $\Delta X(t) : =X(t)-X(t-)= - \Delta Z(t)$ for any $t\ge 0$.
Denote the solution to (\ref{harvest}) with initial condition specified by (\ref{initial}) by $\hat X^{x,\al}\cd$ if necessary.

 We say that $Z$ is an {\em admissible} harvesting strategy if
 \begin{itemize}
 \item[(i)] $Z (t) $ is nonnegative for    any $t\ge 0$ and nondecreasing with respect to $t$,
  %\item[(ii)]   $Z (t)$ is c\`adl\`ag,
 \item[(ii)]    $   \hat X (t )\ge 0$, for any $t\le \tau$, where $\tau$ is the extinction time defined in \eqref{tau} below,
 \item[(iii)] $Z(t)$ is c\`adl\`ag and adapted to ${\cal F}_t:=\sigma\set{w(s),\al(s), 0\le s \le t}$,  and
 \item[(iv)] $J(x,\al,Z) < \infty$ for any $x >0 $ and $\al\in \M $, where $J$ is the functional defined in \eqref{income-J-defn} below.
 \end{itemize}
  Let $\mathcal A$ denote the collection of all admissible harvesting strategies.

Let $f(\cdot, \cdot): \rr_+ \times \M \mapsto \rr_+$ represent the instantaneous marginal
yields accrued from exerting the harvesting strategy $Z$. Assume $f$ is
continuous and non-increasing with respect to $x$. Thus  $f(x,\al) \ge f(y,\al)$ for each $\al\in\M$ whenever $x\le y$.
Moreover, we assume $0< f (0,\al) < \infty$ for
each   $\al \in \M$.  Let $S =(0,\infty)  $, which may be regarded as the {\em survival set} of the species.
Denote the {\em extinction time} by
 \begin{equation}\label{tau} \tau:=\tau^{x,\al}=\inf\set{t\ge 0, \hat X^{x,\al}(t) \notin S}.\end{equation}
  Then for a fixed harvesting process $Z \in \mathcal A$, the expected total discounted value from harvesting is
\begin{equation}\label{income-J-defn}J(x,\al,Z) :=\ex_{x,\al} \int_0^{\tau} e^{-r s}   f(\hat X(s-),\al(s-))  dZ(s)
=\ex  \int_0^{\tau} e^{-r s}   f(\hat X^{x,\al}(s-),\al(s-))  dZ(s),\end{equation}
where  $r \ge 0$ is the discounting factor and   $\ex_{x,\al}$ denotes the expectation with respect to the probability law when the process (\ref{harvest}) starts with initial condition $(x,\al)$ as specified in (\ref{initial}).  The goal is to maximize the expected total discounted value from harvesting and find an optimal harvesting strategy $Z^*$:
\begin{equation}\label{value}
V(x,\al) = J(x,\al,Z^*) := \sup_{Z\in \mathcal{A}} J(x,\al,Z).
\end{equation}

The   dynamic  programming principle takes the form (see \cite{FlemingS,Pha09}):
%\begin{equation}\label{dyna-prog}\begin{aligned}
%V(x,\al)=& \sup_{Z \in \cal A} \ex_{x,\al}\bigg[ \int_0^{\tau\wedge
%  \eta} e^{-rs} f(\hat X(s-),\al(s-)) dZ(s)
%     +     e^{-r(\tau\wedge \eta)}V(\hat X(\tau\wedge\eta),\al(\tau\wedge\eta))\bigg]
%\end{aligned}\end{equation}
\begin{eqnarray}\label{dyna-prog} \nonumber
V(x,\al)&=& \sup_{Z \in \cal A} \ex_{x,\al}\left[ \int_0^{\tau\wedge
  \eta} e^{-rs} f(\hat X^{x,\al}(s-),\al(s-)) dZ(s) \right. \\
& & \qquad \qquad \qquad \left. \rule{0pt}{18pt} +     e^{-r(\tau\wedge \eta)}V(\hat X^{x,\al}(\tau\wedge\eta),\al(\tau\wedge\eta))\right]
\end{eqnarray}
for every $(x,\al)\in S\times\M$ and any stopping time $\eta$.

For later convenience, we introduce the generator of the paired process $(X^{x,\al},\al)$, in which $X^{x,\al}$ satisfies (\ref{dyna-compo}). For any $h(\cdot, \al)\in C^2$, $\al \in\M$, we define
\bed
\op h(x,\al) =   b(x,\al) h'(x,\al) + \frac{1}{2}\sigma^2(x,\al) h''(x,\al) + \sum_{j\in\M} q_{\al j}[h(x,j)-h(x,\al)],
 \eed
where $h'$ and $h''$ denote the first and second order derivatives of $h$ with respect to $x$, respectively.

\section{Verification Theorem and an $\e$-Optimal Policy}\label{sect-verification}

This section establishes a verification theorem whose proof utilizes the generalized It\^o  formula, the monotonicity of $f$, and the regularity of the process $\hat X^{x,\al}$.  We further construct an $\e$-optimal harvesting strategy explicitly in Corollary \ref{cor-1d-optimal} based on the verification theorem and the imposition of additional conditions.

\begin{thm}\label{thm-verification}
 Suppose there exists a function $\phi: S \times \M \mapsto \rr_+$ such that $\phi(\cdot, \al) \in C^2(S) $ for each $\al\in\M$
 and that $\phi$ solves the following coupled system of quasi-variational inequalities:
 \begin{equation}\label{quasi-vari}
\max\set{(\op-r) \phi(x,\al),  f (x,\al)-  \phi'(x,\al)}  = 0,\ \  (x,\al) \in S \times \M,
\end{equation}
where $(\op-r) \phi(x,\al)= \op \phi(x,\al) -r \phi(x,\al)$.
%\bed\begin{aligned} (\op-r) \phi(x,\al) = & \phi'(x,\al)  b(x,\al)  + {1 \over 2}  \phi''(x,\al)\sigma^2(x,\al) % \\ & \ \
%+ \sum_{j\in \M} q_{ij} [\phi(x,j)-\phi(x,\al)] - r \phi(x,\al).\end{aligned}\eed

\begin{itemize}\item[{\em (a)}] Then $\phi(x,\al) \ge V(x,\al)$ for  every $(x,\al)\in S \times \M$.

\item[{\em (b)}] Define the {\em continuation region}
$$ \mathcal C= \set{(x,\al)\in S\times \M:
  f (x,\al)- \phi'(x,\al) < 0}. $$
Assume
 % the following conditions are satisfied:
 % \begin{itemize}\item[\em(i)] $\phi $
 % satisfies $(\op-r) \phi(x,\al)=0$
 % for all $(x,\al)\in \mathcal C$.
 % \item[\em (ii)] Moreover,
there exists a harvesting strategy $ \wdt Z \in \mathcal A$ and corresponding process $\wdt X$ satisfying \eqref{harvest} such that,
 \begin{align}\label{cond-op0} & (\wdt X(t),\al(t)) \in \mathcal C \text { for Lebesgue almost all }  0 \le t \le \tau,   \\
 \label{cond-op1}
& \int_0^t \left[\phi'(\wdt X(s),\al(s)) -f(\wdt X(s),\al(s))\right] d \wdt Z^c(s)  =0, \text{ for any } t \le \tau, \\
\label{cond-op3}&  \lim_{N\to \infty} \ex_{x,\al} \left[e^{-r(\tau\wedge N \wedge \beta_N)} \phi(\wdt X(\tau\wedge N \wedge \beta_N),\al(\tau\wedge N \wedge \beta_N))\right]=0,
\end{align}
and if $\wdt X(s) \not= \wdt X(s-)$, then
\begin{equation}\label{cond-op2}
  \phi(\wdt X(s),\al(s-))-\phi(\wdt X(s-),\al(s-)) = -   f(\wdt X(s-),\al(s-)) \Delta \wdt Z(s),
\end{equation}
 where   $\beta_N:=\inf\{t\ge 0: |\wdt X(t)| \ge N\}$.
Then $\phi(x,\al)= V(x,\al)$ for every $(x,\al )\in S \times \M$ and $\wdt Z $ is an optimal harvesting strategy.
%\end{itemize}
\end{itemize}\end{thm}

\begin{proof} (a) Fix some $(x,\al) \in S \times \M$ and $Z \in \mathcal A$ and let $\hat X$ denote the corresponding solution to (\ref{harvest}).  Choose $N $ sufficiently large so that $\abs{x} < N$ and define $\beta_N:=\inf\set{t\ge 0: | \hat X(t)| \ge N}$.
By virtue of   \cite[Section 2.3]{YZ-10}, \begin{equation}\label{beta-N}\beta_N \to \infty \text{ a.s. as }N \to \infty.\end{equation}
Write $T_N:=N \wedge \beta_N \wedge \tau$. Then It\^{o}'s formula leads to
\bed \begin{aligned}\ex_{x,\al} & [e^{-r T_N}\phi(\hat X(T_N),\al(T_N))] -\phi(x,\al)   \\
= & \ \ex_{x,\al} \int_0^{T_N} e^{-rs} (\op -r) \phi(\hat X(s),\al(s))ds - \ex_{x,\al} \int_0^{T_N} e^{-rs}   \phi' (\hat X(s),\al(s)) dZ^c(s)   \\
 & + \ex_{x,\al} \sum_{0 \le s \le T_N} e^{-rs}\left[\phi(\hat X(s),\al(s-))- \phi(\hat X(s-),\al(s-))
 \right].
\end{aligned}\eed
It follows from (\ref{quasi-vari})  that
\bed \begin{aligned}\ex_{x,\al} & [e^{-r T_N}\phi(\hat X(T_N),\al(T_N))] -\phi(x,\al)
\\ \le  & - \ex_{x,\al} \int_0^{T_N} e^{-rs}   \phi' (\hat X(s),\al(s)) dZ^c(s)
+ \ex_{x,\al} \sum_{0 \le s \le T_N} e^{-rs} \Delta\phi(\hat X(s),\al(s-)),
 %\\  \hfill & \hspace{11em}
 % + \sum_{i=1}^n  \frac{\partial \phi}{\partial x_i}(\hat X(s-),\al(s-))
 % \Delta Z_i(s)) \Big],
\end{aligned}\eed
where $\Delta\phi(\hat X(s),\al(s-))= \phi(\hat X(s),\al(s-))- \phi(\hat X(s-),\al(s-)) $. Apply the mean value theorem to $\Delta\phi(\hat X(s),\al(s-))$ and we obtain
$$\Delta\phi(\hat X(s),\al(s-)) =  \phi'(\xi(s),\al(s-)) \Delta \hat X(s)= -  \phi'(\xi(s),\al(s-)) \Delta Z(s),$$ where $\xi(s)= \theta(s) \hat X(s) + (1-\theta(s)) \hat X(s-)$
for some $\theta(s) \in (0,1)$. Note that  $\hat X(s) \le \xi(s) \le \hat X(s-)$.
Thus it follows that
\bed \begin{aligned}\phi(x,\al) \ge &  \ex_{x,\al}  [e^{-r T_N}\phi(\hat X(T_N),\al(T_N))]+  \ex_{x,\al} \int_0^{T_N} e^{-rs} \phi'(\hat X(s),\al(s)) dZ^c(s)
\\ &  +  \ex_{x,\al} \sum_{0 \le s \le T_N} e^{-rs}    \phi'(\xi(s),\al(s-)) \Delta Z(s).  \end{aligned}\eed
Using (\ref{quasi-vari}) again
and noting that $\phi$ is nonnegative and that $f(\cdot,\al)$ is nonincreasing for each $\al \in \M$, it follows that
\bed \begin{aligned}
\phi(x,\al) & \ge \ex_{x,\al} \int_0^{T_N} e^{-rs}   f (\hat X(s),\al(s)) dZ^c(s) +  \ex_{x,\al} \sum_{0 \le s \le T_N} e^{-rs}   f(\xi(s),\al(s-)) \Delta Z(s) \\
& \ge \ex_{x,\al} \int_0^{T_N} e^{-rs}   f(\hat X(s),\al(s)) dZ^c(s)% \\ & \hfill \hspace{4em}
+  \ex_{x,\al} \sum_{0 \le s \le T_N} e^{-rs}   f(\hat X(s-),\al(s-)) \Delta Z(s)\\
&= \ex_{x,\al}\int_0^{T_N} e^{-rs} f(\hat X(s-),\al(s-)) dZ (s).
\end{aligned}\eed
Now letting $N \to \infty$, it follows from (\ref{beta-N}) and the bounded convergence theorem that
$$ \phi(x,\al) \ge \ex_{x,\al} \int_0^\tau e^{-rs}   f(\hat X(s-),\al(s-))  dZ(s). $$
Finally, taking supremum over   all $Z \in \cal A$, we obtain $\phi(x,\al) \ge V(x,\al)$, as desired.

(b) Let $\wdt Z \in \cal A$ satisfy (\ref{cond-op0})--(\ref{cond-op2}). Define $\beta_N$ and $T_N$ as before with $\wdt X$ replacing $\hat X$. As in part (a), we have from It\^o's formula that
\bed \begin{aligned}\ex_{x,\al} & [e^{-r T_N}\phi(\wdt X(T_N),\al(T_N))] -\phi(x,\al)   \\
= & \ \ex_{x,\al} \int_0^{T_N} e^{-rs} (\op -r) \phi(\wdt X(s),\al(s))ds - \ex_{x,\al} \int_0^{T_N} e^{-rs}  \phi'(\wdt X(s),\al(s)) d\wdt Z^c(s)   \\
 & + \ex_{x,\al} \sum_{0 \le s \le T_N} e^{-rs}\left[\phi(\wdt X(s),\al(s-))- \phi(\wdt X(s-),\al(s-))
 \right].
\end{aligned}\eed
By (\ref{cond-op0}), $(\op-r) \phi(\wdt X(s),\al(s)) =0$ for almost all $s \in [0,\tau]$. This, together with  (\ref{cond-op1}) and (\ref{cond-op2}), implies that
\bed \phi(x,\al) = \ex_{x,\al}   [e^{-r T_N}\phi(\wdt X(T_N),\al(T_N))]
+  \ex_{x,\al}\int_0^{T_N} e^{-rs}   f(\wdt X(s-),\al(s-)) d\wdt Z(s). \eed
Letting $N\to \infty$ and using (\ref{cond-op3}) and (\ref{beta-N}),
 we obtain
$$\phi(x,\al)= \ex_{x,\al}\int_0^{\tau} e^{-rs}   f(\wdt X(s-),\al(s-))  d\wdt Z(s) . $$ This shows that $\phi(x,\al)= V(x,\al)$ for every $(x,\al )\in S \times \M$ and $\wdt Z$ is an optimal harvesting strategy.
\end{proof}

\begin{rem}
The conditions of \thmref{thm-verification} can be weakened. In fact, by virtue of
  %Fleming and Soner \cite{FlemingS} and {\O}ksendal
  \cite[Appendix D]{Oksendal}, we need only to assume that (i) $\phi(\cdot,\al) \in C^1(S)\cap C^2(S-D)$ for each $\al\in \M$, where $D$ is countable set of points, and (ii) $  \phi''(x+) < \infty$, $\phi''(x-) < \infty$ for all $x\in D$.
Under these conditions, there exist sequences $\set{\phi_j(\cdot, \al)}_{j=1}^\infty$, $\al \in \M$ such that $\phi_j(\cdot, \al)\in
C^2(S) $ for each $\al\in\M$. Moreover, the following are satisfied:
\begin{itemize}
\item[{(a)}] for each $\al \in \M$,  $\lim_{j\to \infty}\phi_j(\cdot, \al) \to \phi(\cdot, \al)$ uniformly on compact subsets of $S$,
\item[(b)]  $\lim_{j\to \infty}(\op -r )\phi_j(x, \al) \to \phi(x, \al)$ uniformly on compact subsets of $S-D$, $\al \in \M$, and
\item[(c)] $\set{(\op-r)\phi_j}_{j=1}^\infty$ is locally bounded on $S\times \M$.
\end{itemize}
Then, we can first work with the sequence $\phi_j$ exactly the same way as in the proof of \thmref{thm-verification}. Next, using (a), (b), and (c), we can pass to the limit as $j\to \infty$ to obtain the same conclusions. The reader is referred to   \cite{Oksendal} for details.
\end{rem}

By virtue of  \thmref{thm-verification}(a),  any sufficiently smooth solution to
% the coupled system of quasi-variational inequalities
\eqref{quasi-vari} is an upper bound for the value function $V$.
Further, the additional conditions in \thmref{thm-verification}(b) will help us to  find an optimal harvesting strategy. In practice, it is, however, usually very hard to find an explicit solution to
\eqref{quasi-vari}. In particular, with the presence of regime switching, \eqref{quasi-vari} is a coupled system of quasi-variational inequalities,
  a closed form solution is virtually impossible except in some special cases (see Examples \ref{example-2regimes} and \ref{example-non-constant-price}).
Nevertheless, some results about the value function can be derived; \propref{prop-1d} below gives an upper bound for $V$ when $f(\cdot,\al)$ is smooth for each $\al$. Furthermore, under additional assumptions, we explicitly construct an $\e$-optimal harvesting strategy in \thmref{cor-1d-optimal}.

For any $x>0$ and $\al\in \M$, define
\begin{equation}\label{g-defn} g (x,\al)= \int_0^x f(y,\al)dy.\end{equation}
Then it follows that $g  $ is nonnegative and $g '(x,\al)= f(x,\al)$. Moreover, if $f(\cdot,\al) \in C^1(S)$, then $g ''(x,\al)= f'(x,\al)\le 0$ because $f(\cdot, \al)$ is nonincreasing for each $\al\in \M$. This shows that $g (\cdot, \al)$ is concave for each $\al \in \M$.
\begin{prop}\label{prop-1d} Assume that $f(\cdot, \al)\in C^1(S)$ and $f(\cdot, \al)$ is non-increasing for each $\al\in\M$.
     Then we have
\begin{equation}\label{value-ineq1} V(x,\al) \le g (x,\al) + \sup_{Z\in \cal A} \ex_{x,\al} \int_0^\tau e^{-rs} (\op-r) g (\hat X(s),\al(s))ds. \end{equation}
\end{prop}

\begin{proof} Fix some $(x,\al) \in S \times \M$ and $Z \in \mathcal A$ and let $\hat X$ denote the corresponding solution to (\ref{harvest}).  Let $T_N$ be as in the proof of \thmref{thm-verification}. Apply It\^o's formula using $g$ to obtain
\bed\begin{aligned} & \ex_{x,\al} \left[e^{-rT_N}  g (\hat X(T_N),\al(T_N))\right] - g (x,\al)\\ & \ \ = \ex_{x,\al} \int_0^{T_N} e^{-rs} (\op-r) g (\hat X(s),\al(s)) ds
 - \ex_{x,\al}  \int_0^{T_N} e^{-rs} g '(\hat X(s-),\al(s-))dZ(s)
 \\ &  \ \ \quad + \ex_{x,\al} \sum_{0\le s \le T_N}  e^{-rs} \Big[ g (\hat X(s),\al(s-))- g (\hat X(s-),\al(s-)) - g '(\hat X(s-),\al(s-)) \Delta \hat X(s)\big].
 \end{aligned}\eed
Since $g (\cdot, \al)$ is concave for each $\al \in \M$, it follows that $$g (\hat X(s),\al(s-)) \le  g (\hat X(s-),\al(s-)) +  g '(\hat X(s-),\al(s-)) (\hat X(s)-\hat X(s-)).$$
 Thus we have
 \bed g (x,\al) \ge  \ex_{x,\al}  \int_0^{T_N} e^{-rs} f(\hat X(s-),\al(s-))dZ(s) - \ex_{x,\al} \int_0^{T_N} e^{-rs} (\op-r) g (\hat X(s),\al(s)) ds .\eed
Now letting $N\to \infty$ and using the same argument as in the proof of \thmref{thm-verification}, we obtain
\bed g (x,\al) \ge  \ex_{x,\al}  \int_0^{\tau} e^{-rs} f(\hat X(s-),\al(s-))dZ(s) - \ex_{x,\al} \int_0^{\tau} e^{-rs} (\op-r) g (\hat X(s),\al(s)) ds ,\eed
from which (\ref{value-ineq1}) follows by taking supremum over $Z \in \cal A$.
\end{proof}

\begin{thm}\label{cor-1d-optimal} Assume, in addition to the conditions of \propref{prop-1d}, $(\op-r) g(x,\al) \le 0$ for all $(x,\al)\in S\times \M$.
\begin{itemize}\item[{\em (i)}]  Suppose that there exists a constant $L>0$ such that
\beq{f-lip} \abs{f(x,\al)-f(y,\al)} \le L \abs{x-y}, \hbox{ for all }x,y\in S \hbox{ and } \al\in \M. \eeq  Then for any $\e > 0$, there exists a
harvesting strategy $Z^\e\in \cal A$
under which \begin{equation}\label{e-optimal-to-g}g(x,\al) -\e \le J(x,\al,Z^\e) \le V(x,\al) \le g(x,\al).
\end{equation} The harvesting strategy $Z^\e $ is a ``chattering policy'' that instantaneously harvests a sufficiently small amount many times in a sufficiently small interval of time until the species becomes extinct.
\item[{\em (ii)}]  In particular, if $f(x,\al)\equiv  f(\al)$ for all $x\in S$ and each $\al\in \M$, then
   \begin{equation}\label{optimal-to-g}
   V(x,\al)= g(x,\al)= f(\al) x,
   \end{equation}
  and the optimal harvesting strategy  is to drive the process instantaneously to extinction.
\end{itemize}
\end{thm}

\begin{proof}
Fix some $(x,\al)\in S \times \M$.
By \propref{prop-1d} and the condition $(\op-r) g(x,\al) \le 0$, we have \begin{equation}\label{V-ge-g}V(x,\al) \le g(x,\al).\end{equation}
The rest of proof is divided into two parts. Part 1 is devoted to the proof of \eqref{e-optimal-to-g} while the second part establishes \eqref{optimal-to-g}.

\noindent
{\em Part 1.}\/ Since $f(\cdot,\al)$ is continuous, for any $\e >0$,  there exists an $N \in \mathbb N$ such that   $$  R(f,\al):=\sum_{i=0}^{n-1} f(x_i,\al) \delta > \int_0^x f(y,\al)dy -\e/3 =g(x,\al)-\e/3, \hbox{ for any } n \ge N, $$ where $\delta= x/n$ and $x_i= x- i\delta$, $i=0,1,\dots, n-1$.
 Note that we also have $  R(f,\al)\le  g(x,\al)$ since $f(\cdot,\al)$ is non-increasing. Thus it follows that
 \begin{equation}\label{R-g-difference} \abs{R(f,\al) -g(x,\al)} < \e/3. \end{equation}

 Let $\varsigma  = n^{-5}$ and $t_i = i \varsigma /n$, $i=0,1, \dots, n $.  We construct a harvesting strategy $Z= Z^\e$ which increases only on the set $\{t_i:i=0,\ldots,n\}$; denote the corresponding harvested process by $\hat X$. Note that $\hat X(t_0-) = x_0=x$. Define $\Delta Z(t_0) =Z(t_0) = \delta$, observe
 $\hat X(t_0) = x_1$ and it therefore follows that
 $$\hat X(t_1-) = X(t_0) + \int_{t_0}^{t_1} b(\hat X(s),\al(s))ds + \int_{t_0}^{t_1} \sigma(\hat X(s),\al(s))dw(s).$$
At time $t=t_1$, define $\Delta Z(t_1) = (\hat X(t_1-)-x_2)^+$ so that $\hat X(t_1) \leq x_2$ and allow the process $\hat X$ to diffuse until time $t=t_2$.
 In general, for $i=1, \dots, n-1$, we define
 $$ \Delta Z(t_i) = \(\hat X(t_i-)- x_{i+1}\)^+$$
so that
$$\hat X(t_i)=\hat X(t_i-) -\Delta Z(t_i),$$
and
$$ \hat X(t_{i+1}-)= \hat X(t_i) + \int_{t_i}^{t_{i+1}} b(\hat X(s),\al(s))ds +  \int_{t_i}^{t_{i+1}}  \sigma(\hat X(s),\al(s))dw(s).  $$
 Note that $\hat X(t_i) = x_{i+1}$ if $\Delta Z(t_i) > 0$.
 The expected total discounted income from the harvesting strategy $Z$ is $$J(x,\al, Z) = \ex_{x,\al} \sum_{i=0}^{n-1} e^{-r t_i} f(\hat X(t_i-),\al(t_i-)) \Delta Z(t_i).$$ Next we want to show that $\abs{J(x,\al,Z)- R(f,\al)} < \e /3$.
 In fact, we have
 \begin{equation}\label{JR-diff} \barray\abs{J(x,\al,Z)- R(f,\al)} \ad \le \sum^{n-1}_{i=0} \ex_{x,\al}\abs{e^{-rt_i} f(\hat X(t_i-),\al(t_i-))\Delta Z(t_i) -f(x_i,\al)\delta}\\
 \ad  \le \sum^{n-1}_{i=0}\bigg[ \ex_{x,\al} \abs{[f(\hat X(t_i-),\al(t_i-))-f(x_i,\al)]\delta} \\ \aad \ \qquad  + \ \ex_{x,\al}\abs{f(\hat X(t_i-),\al(t_i-))[\Delta Z(t_i)-\delta]} \\ \aad \ \qquad +\  \ex_{x,\al}\abs{[e^{-r t_i} -1]f(\hat X(t_i-),\al(t_i-))\Delta Z(t_i)} \bigg] \\
 \ad =: \sum_{i=1}^n \(A_i + B_i + C_i\). \earray \end{equation}
In the following, we analyze the terms $A_i$, $B_i$, and $C_i$ separately. To this end,  for any $i =0, 1, \dots, n-1$, we apply \cite[Proposition 2.3]{YZ-10} to obtain
\begin{align}\label{moment-estimate}
& \ex \abs{\int_{t_i}^{t_{i+1}} b(\hat X(s),\al(s))ds + \int_{t_i}^{t_{i+1}} \sigma(\hat X(s),\al(s))dw(s) }^2 \le K (t_{i+1}-t_i) = Kt_1,
\\ \label{DZ-estimate} & \ex \abs{\Delta Z(t_i)} =\ex \abs{(\hat X(t_i-)-x_{i+1})^+} \le K,
\end{align} where $K$ is a generic  positive constant  depending only on $x$, $m$, and the constant $\kappa_0$ in \eqref{ito-condition}.
Also, in the sequel, the exact value of $K$ may change in different appearances.    Then it follows from the Tchebychev inequality that
\beq{ineq-prob-Z=0} \barray
\pr \set{\Delta Z(t_1) =0} \ad = \pr \set{\hat X(t_1-) \le x_2} \\ \ad
= \pr \set{\int_{t_0}^{t_1} b(\hat X(s),\al(s))ds + \int_{t_0}^{t_1} \sigma(\hat X(s),\al(s))dw(s) \le -\delta }\\
\ad \le \pr \set{\abs{\int_{t_0}^{t_1} b(\hat X(s),\al(s))ds + \int_{t_0}^{t_1} \sigma(\hat X(s),\al(s))dw(s) } \ge \delta } \\
\ad \le \frac{Kt_1}{\delta^2}. \earray \eeq
 Note that $\hat X(t_1) =x_2$ if $\Delta Z(t_1) > 0$. Thus we have
 \beq{ineq-prob-X-1} \pr \set{\hat X(t_1) \not = x_2} \le \pr \set{\Delta Z(t_1) =0} \le \frac{Kt_1}{\delta^2}.\eeq
 Using the same arguments as those in \eqref{ineq-prob-Z=0} and \eqref{ineq-prob-X-1}, we have \begin{align} \label{ineq-prob-Z-2}    \barray \pr\set{\Delta Z(t_2)=0}\ad = \pr\set{\Delta Z(t_2) =0,\hat X(t_1) =x_2}+ \pr \set{\Delta Z(t_2) =0 ,\hat X(t_1)\not= x_2}\\
 \ad \le  \frac{Kt_1}{\delta^2} + \frac{Kt_1}{\delta^2} = \frac{Kt_2}{\delta^2},\earray\end{align}
 and \begin{align}
 \label{ineq-prob-X-2}   \pr\set{\hat X(t_2) \not= x_3} \le \pr\set{\Delta Z(t_2)=0} \le \frac{Kt_2}{\delta^2}.  \end{align}
 Continuing in this manner, it follows that for any $i=1, 2, \dots, n-1$,
 \begin{align} \label{ineq-prob-Z-i}
 \pr\set{\Delta Z(t_i)=0} \le \frac{Kt_i}{ \delta^2}, \quad \mbox{ and }\\
 \label{ineq-prob-X-i} \pr\set{\hat X(t_i) \not= x_{i+1}} \le \frac{Kt_i}{ \delta^2}.
 \end{align}
 Using the conditions that $f$ is Lipschitz continuous and uniformly bounded,  we compute
 \bed \barray
  A_i \ad \le  \ex \abs{f(\hat X(t_i-), \al) -f(x_i,\al)} \delta + \ex \abs{f(\hat X(t_i-), \al(t_i-))- f(\hat X(t_i-), \al) }\delta \\
  \ad \le L \ex \abs{\hat X(t_i-) -x_i} \delta + K \pr\set{\al(t_i-)\not =\al} \delta \\
  \ad \le L \ex \abs{\hat X(t_{i}-) -x_i} \delta  + K t_i \delta, \earray \eed
  where in the last inequality, we used \eqref{Q-gen}.
  But using \eqref{ineq-prob-X-i}, \eqref{moment-estimate}, and \cite[Proposition 2.3]{YZ-10}, we obtain
  \bed \barray \ex\abs{\hat X(t_i-) -x_i} \ad \le \ex \abs{\hat X(t_{i-1}) -x_i} + \ex\abs{\int_{t_{i-1}}^{t_i}b(\hat X(s),\al(s))ds + \int_{t_{i-1}}^{t_i}\sigma(\hat X(s),\al(s))dw(s) } \\
  \ad \le \ex^{1/2} \abs{\hat X(t_{i-1}) -x_i}^2\ex^{1/2}[I_{\set{\hat X(t_{i-1}) \not= x_i}}] + Kt_1 \\
  \ad \le \frac{K\sqrt{t_{i-1}}}{ \delta} + Kt_1 \le  \frac{K\sqrt{t_{i}}}{ \delta} + Kt_1. \earray \eed
  Thus it follows that
  \begin{equation}\label{ineq-A-estimate}
A_i \le K \delta\(\frac{ \sqrt{ t_i}}{ \delta} + t_i +  t_1 \)= K (\sqrt{t_i}+ t_i  \delta+ t_1 \delta).
  \end{equation}

Next we estimate $B_i$. Since $f$ is uniformly bounded, it follows that
\bed \barray B_i \ad \le K \ex\abs{\Delta Z(t_i) -\delta}\\
 \ad = K \ex \abs{(\Delta Z(t_i) -\delta)I_{\set{\Delta Z(t_i)=0}}} + K \ex \abs{(\Delta Z(t_i) -\delta)I_{\set{\Delta Z(t_i)\not=0}}I_{\set{\hat X(t_{i-1})=  x_i}}} \\ \aad \quad +  K \ex \abs{(\Delta Z(t_i) -\delta)I_{\set{\Delta Z(t_i)\not=0}}I_{\set{\hat X(t_{i-1})\not= x_i}}}
 \\ \ad := B_{i1}+B_{i2}+ B_{i3}. \earray \eed
   Note that \eqref{ineq-prob-Z-i} implies that $B_{i1} \le \delta \frac{Kt_i}{ \delta^2}= \frac{Kt_i}{ \delta}.$
   Using the definition of $\Delta Z(t_i)$ and \eqref{moment-estimate}, we have
   \bed \barray B_{i2} \ad = K \ex \abs{(\hat X(t_i-)-x_{i+1} -\delta)
   I_{\set{\Delta Z(t_i)\not=0}}I_{\set{\hat X(t_{i-1})=  x_i}}}\\
   \ad = K \ex \abs{\( %\hat X(t_{i-1})  -x_i +
    \int_{t_{i-1}}^{t_i} b(\hat X(s),\al(s))ds +\int_{t_{i-1}}^{t_i} \sigma(\hat X(s),\al(s))dw(s) \)I_{\set{\Delta Z(t_i)\not=0}}I_{\set{\hat X(t_{i-1})=  x_i}} }\\
   \ad \le  Kt_1. \earray \eed
 Concerning the term $B_{i3}$, we use the Cauchy-Schwartz inequality,
  \eqref{DZ-estimate}, and \eqref{ineq-prob-X-i}:
 \bed \barray B_{i3} \ad \le K \ex^{1/2} \abs{(\Delta Z(t_i)-\delta)I_{\set{\Delta Z(t_i) \not=0}}}^2 \ex^{1/2} \abs{I_{\set{\hat X(t_{i-1})\not= x_i}}}^2 \\
 \ad \le  K \frac{\sqrt{t_{i-1}}}{\delta} \le K \frac{\sqrt{t_{i}}}{\delta}.
\earray \eed
Putting these estimates together, we obtain
\begin{equation}\label{ineq-B-estimate}
B_i \le K \frac{t_i}{\delta} + Kt_1 + K  \frac{\sqrt{t_i}}{\delta}\le K(t_1 +\frac{\sqrt{t_i}}{\delta} ).
\end{equation}

  For the  term $C_i$, we again use the   uniform boundedness of $f$ and
  %\cite[Proposition 2.3]{YZ-10}
  \eqref{DZ-estimate} to obtain
  \begin{equation}\label{ineq-C-estimate}
  C_i =\ex \abs{(e^{-rt_i}-1) f(\hat X(t_i-),\al(t_i-)) \Delta Z(t_i)} \le K(1-e^{-r t_i}) = Kr t_i + o(t_i) \le Kt_i.\end{equation}

Now using the estimates \eqref{ineq-A-estimate}, \eqref{ineq-B-estimate}, and \eqref{ineq-C-estimate} in \eqref{JR-diff}, and noting $\delta= x n^{-1}$,  $t_i= it_1$, and $t_1= n^{-6}$,
\bed \barray \abs{J(x,\al,Z) - R(f,\al)} \ad \le \sum_{i=1}^n \(K (\sqrt{t_i}+ t_i  \delta+ t_1 \delta)+ K(  t_1 + \frac{\sqrt{t_i}}{\delta}) + Kt_i\)\\
\ad \le K\(t_1\sum^n_{i=1} ( i +1)+ \sqrt{t_1} \frac{\delta +1 }{\delta} \sum^n_{i=1} \sqrt{i} \) \\
\ad \le K \(n^2 t_1   + \sqrt{t_1} \frac{\delta +1 }{\delta} n^{3/2}\)\\
\ad \le K n^2 n^{-6} + n^{-3} n n^{3/2} \le K n^{-1/2}.
 \earray \eed
 Finally we choose $n$ sufficiently large so that \eqref{R-g-difference} holds and $\abs{J(x,\al,Z) - R(f,\al)} < \e/3$. Then it follows that
 \bed \abs{J(x,\al,Z)- g(x,\al)} \le \abs{J(x,\al,Z)- R(f,\al)}+ \abs{R(f,\al)-g(x,\al)} < \e/3 + \e/3 < \e.\eed
Now \eqref{e-optimal-to-g} follows in view  of \eqref{V-ge-g}.

\noindent
{\em Part 2.}\/ If   $f(x,\al)\equiv  f(\al)$ for all $x\in S$ and each $\al\in \M$, then we choose $Z$ to be the harvesting policy which drives the process $\hat X$ instantaneously from state $x$ to state $0$. It follows that $\tau=0$ and
$$ J(x,\al,Z) = f(\al)x = g(x,\al)= V(x,\al).$$
  This finishes the proof.
\end{proof}

\begin{rem}\label{rem-Alvarez}
In \cite[Corollary 1]{Alvarez}, it was commented that ``if the convenience yield from holding reserves is non-positive at all states then the optimal
policy is to deplete the reserves at an infinitely fast rate but only in small proportions at a time (a form of a `chattering policy').'' While the intuition in \cite{Alvarez} is correct,  the optimal
policy in \cite{Alvarez} is not admissible in our context, because it is not well defined at time  $0$. In \thmref{cor-1d-optimal}, we explicitly constructed an admissible   harvesting
policy, under which the expected total discounted income from harvesting is $\e$-optimal.
\end{rem}

%\begin{thm}\label{thm-about-e-optimal}
%Suppose
%\begin{itemize}
%\item[{\em (i)}] there exists a sufficiently smooth function $\phi: S\times \M \mapsto \rr_+ $ that solves
%  %the coupled system of quasi-variational inequalities
%  \eqref{quasi-vari},

% \item[{\em (ii)}] for some $\alpha \in \M$ and some constant $0\le b_\alpha < \infty$,
%$(0,b_\alpha) \subset \mathcal C$ and  $(b_\alpha, \infty) \subset \mathcal H$,

% \item[{\em (iii)}] the function $f(\cdot, \alpha)$ is Lipschitz continuous and strictly decreasing,

%\end{itemize} Then for any $\e >0$, there is an admissible harvesting strategy $Z^\e$ such that
%$$\phi(x,\alpha) -\e < J(x,\alpha, Z^\e) \le \phi(x,\alpha), \text{ for any }x \in (b_\alpha, \infty).$$
% The harvesting strategy $Z^\e$ is a ``chattering policy'' as described in \thmref{cor-1d-optimal},
%which harvest a sufficiently small amount instantaneously for sufficiently large number of times.
%\end{thm}
%\begin{proof}
%For any $x\in (b_\alpha, \infty) $, put $\delta = \frac{x-b_\alpha}{n}$ where $n \in \mathbb N$.
% The rest of proof is similar to the argument in \thmref{cor-1d-optimal}.
%\end{proof}

\section{Examples}\label{sect-example}

We  provide two examples
 to demonstrate our results in the previous section. They reveal that in the setting of regime switching, it is much harder to obtain the value functions and the ($\e$-) optimal harvesting strategies due to the coupling in the system of quasi variational inequalities.

\begin{exm}\label{example-2regimes}
We assume the growth of a certain species (or a certain risky investment) is governed by a regime switching geometric Brownian motion
\beq{exam-gbm} dX(t)= \mu(\al(t)) X(t)dt + \sigma(\al(t))X(t)dw(t), \eeq
and the harvested process is given by
\beq{exam-harv} d\hat X(t)= \mu(\al(t)) \hat X(t)dt + \sigma(\al(t))\hat X(t)dw(t)-dZ(t), \eeq
where $Z(t)$ denotes the total amount  of harvest (or dividends) up to time $t$,  $w$ is a standard Brownian motion,
 $\al$ is a continuous time Markov chain with state space $\M=\set{1,\dots,m}$, and for each $\al \in \M$, $\mu_\al=\mu(\al)$ and $\sigma_\al=\sigma(\al)$ are constants. Our objective is to maximize the expected discounted income from harvest and  find an optimal harvesting policy, i.e., we want to find
 \beq{objec-1d} V(x,\al) = J(x,\al,Z^*) =\sup_{Z\in\cal A}\ex_{x,\al} \int^\tau_0 e^{-rs} dZ(s),\eeq
 where $r>0$ is the discount factor.  Note that $f \equiv 1$ in this example.

First consider the case when $m=1$; that is, there is only a static environment so no regime switching occurs. It is clear that if $\mu > r$, then $V(x,1)=\infty$. On the other hand, if $\mu\le r$, then $V(x,1)=x$ and the optimal harvesting policy is to drive the process instantaneously to extinction ($\tau=0$ a.s.).  We refer the reader to   \cite{Alvarez} or  \cite{A-Shepp} for details.

 Now let $m=2$ and assume that the continuous time Markov chain $\al$ is generated by $Q=\begin{pmatrix} -\lambda_1 & \lambda_1 \\  \lambda_2 & -\lambda_2\end{pmatrix}$, where $\lambda_1>0$ and $\lambda_2 >0$. Without loss of generality, we further assume that $\mu_1 \le \mu_2$.

 Case 1: $\mu_1 \le r$ and $\mu_2 \le r$. In this case, we have $g(x,1)=g(x,2)=x$ and
 $$(\op-r) g(x,i) = \mu_i x -(\lambda_i+ r ) x + \lambda_i x = (\mu_i -r) x \le 0, \ \ x>0, \ i=1,2.$$
 Then \thmref{cor-1d-optimal} implies that $V(x,1)=V(x,2)=x$ and that the optimal policy is to drive the process instantaneously to 0.

 Case 2: $\mu_1  < r < \mu_2 \le  \xi $, where
 \begin{equation}\label{eq:exm1-xi-defn}\xi= \frac{r\lambda_1 + (r-\mu_1)(r+ \lambda_2)}{r+ \lambda_1 -\mu_1} .\end{equation} Note that $\xi > r$.
   %and $(r-\mu_1) (\lambda_2 + r -\mu_2) \ge (\mu_2-r)\lambda_1$.
In this case, it can be shown that the unique solution to the system of coupled quasi variational inequalities
\bed \max\set{(\op-r) \phi(x,\al), 1- \phi'(x,\al)}=0, \ \  x> 0, \ \al=1,2,\eed
is \beq{soln-1d-ex} \phi(x,1)=x,\ \ \phi(x,2)=\frac{\lambda_2}{\lambda_2+r-\mu_2}x. \eeq Therefore \thmref{thm-verification} implies  that  $\phi(x,\al) \ge V(x,\al)$, $\al=1,2$.

Next we show that there is a harvesting strategy $Z^*$ under which $J(x,\al,Z^*)=\phi(x,\al)$  for all $ x>0$ and $\al=1,2$.
To this end, we denote  the  harvesting region  by $\mathcal H:= (0,\infty) \times \set{1}$ and the continuation region by  $\mathcal C:= (0,\infty) \times \set{2}$.
Let  the harvesting policy $Z^*$ be such that %by The optimal policy is to
it  drives the process instantaneously to the origin once the Markov chain enters state 1 or the process enters the harvesting region. Consequently, the extinction time is $$\tau=\tau^{x,\al}=\inf\set{t\ge 0: \al(t) =1}=\inf\set{t\ge 0: (\hat X^{x,\al}(t),\al (t))\in \mathcal H}.$$
One can verify that this harvesting policy and the corresponding harvested process satisfy all conditions in \thmref{thm-verification}, part (b).
 In fact, it is obvious   that $J(x,1,Z^*)=x=\phi(x,1)$. Next we consider $J(x,2,Z^*)$. Note that $\tau = \tau^{x,2}$ has exponential distribution with parameter $\lambda_2$ and that $$ \hat X (t) =\hat X^{x,2} (t)  = x \exp\set{(\mu_2 -\frac{1}{2}\sigma_2^2)t + \sigma_2 w(t)}, \text{ for all } t\in [0,\tau].$$ Therefore it follows that
 \bed \barray J(x,2,Z^*)\ad = \ex_{x,2} \int_0^{\tau} e^{-rs} dZ^*(s) = \ex_{x,2}  [e^{-r\tau}\hat X(\tau)] \\ \ad =\int_0^\infty e^{-rt} x \exp\set{ (\mu_2-\frac{1}{2} \sigma_2^2)t + \frac{1}{2}\sigma_2^2 t} \lambda_2 e^{-\lambda_2 t}dt \\\ad =\frac{\lambda_2}{\lambda_2+r-\mu_2}x = \phi(x,2). \earray \eed
Hence $V(x,\al)=\phi(x,\al)$ for all $(x,\al) \in (0,\infty) \times \set{1,2}$ and $Z^*$ is an optimal harvesting strategy.

 Case 3: $\mu_1 < r < \xi < \mu_2$, where $\xi$ is defined in \eqref{eq:exm1-xi-defn}. We claim that $V(x,1)=V(x,2) =\infty$ for any $x>0$.
 This is quite interesting. It indicates that even though $\mu_1 < r$, we still have $V(x,1) =\infty$ thanks to the  switching component $\al$.
 In the sequel, we demonstrate that there exists an admissible harvesting policy $Z$
  under which $J(x,1,Z)$ and $J(x,2,Z)$ can be arbitrarily large and hence the claim follows.

Fix some $M>0$ and define $\eta:= \inf\set{t\ge 0: (X(t),\al(t))=(M,2)}$. Note that the function $u(x,i)= \ex_{x,i}[e^{-r\eta}]$, $0<x< M$, $i=1,2$, solves the differential equation $(\op-r) u(x,i)=0$. That is, $u$ is a solution to the coupled system of differential equations
\begin{equation}\label{exm-2regimes} \barray
\ad \frac{1}{2} \sigma_1^2 x^2 u''(x,1) + \mu_1 x u'(x,1) -(r+\lambda_1 ) u(x,1) + \lambda_1 u(x,2)=0,\\[1.5ex]
\ad \frac{1}{2} \sigma_2^2 x^2 u''(x,2) + \mu_2 x u'(x,2) -(r+\lambda_2 ) u(x,2) + \lambda_2 u(x,1)=0.\earray\end{equation}
The characteristic equation of (\ref{exm-2regimes}) is $h(x)=g_1(x)g_2(x)-\lambda_1 \lambda_2=0$, where $g_i(x)= \frac{1}{2} \sigma_i^2 x(x-1) + \mu_i x -r-\lambda_i$, $i=1,2$. As argued in   \cite{Guo-Z}, $h(x)$ has four real roots $\beta_1 > \beta_2 >0 > \beta_3 > \beta_4$. Moreover, the condition
$\mu_1 < r < \xi  <  \mu_2$ implies that
$ 1> \beta_2 >0$.
Therefore $u(x,i)$, a solution of (\ref{exm-2regimes}), can be written as
$$u(x,i)= \ex_{x,i}[e^{-r\eta}]= \sum_{j=1}^4 C^i_j x^{\beta_j}, \ \ 0< x < M, \ i=1,2$$ for some constants $C_j^i$, $i=1,2$ and $j=1,2,3,4$. As noted in \cite{Guo-Z},   $C_j^2=l_j C^1_j$, where $l_j =-\frac{\lambda_2}{g_2(\beta_j)}= -\frac{g_1(\beta_j)}{\lambda_1}$.
But as $$X(t) =x \exp\set{\int_0^t [\mu(\al(s))-\frac{1}{2} \sigma^2(\al(s))]ds + \int_0^t \sg(\al(s))dw(s)},$$
it follows that $\eta \to \infty $ a.s. as $x \downarrow 0$. Thus   for $i=1,2$, we have $\ex_{x,i} [e^{-r\eta}] \to 0$ as $x\downarrow 0$
and hence
$$\barray \ex_{x,1} [e^{-r\eta}]= C_1 x^{\beta_1} + C_2 x^{\beta_2}, \ \ %\\
  \ex_{x,2} [e^{-r\eta}]= l_1 C_1 x^{\beta_1} + l_2 C_2 x^{\beta_2}, \earray $$
where $C_1$ and $C_2$ are constants, and $$%\barray
l_1=- \dfrac{\lambda_2}{g_2(\beta_1)}=-\dfrac{g_1(\beta_1)}{\lambda_1} < 0,\ \ l_2=- \dfrac{\lambda_2}{g_2(\beta_2)}=-\dfrac{g_1(\beta_2)}{\lambda_1} > 0.$$
Now the boundary conditions yield
\bed \barray \ex_{M,1} [e^{-r \eta}] =c = C_1 M^{\beta_1} + C_2 M^{\beta_2}, \ \
 \ex_{M,2}[e^{-r \eta} ]=1 =l_1 C_1 M^{\beta_1} + l_2 C_2 M^{\beta_2},\earray \eed
 where $0< c\le 1$. Solve the above equations for $C_1$ and $C_2$ and we obtain
 $$C_1= \frac{l_2 c-1}{(l_2-l_1)M^{\beta_1}}, \ \
   C_2= \frac{1- l_1 c}{(l_2-l_1)M^{\beta_2}}.$$
 Notice that $C_2 > 0$. Consequently, we can write for $x\in (0,M)$ that
 $$\ex_{x,1}[e^{-r \eta} ] = \frac{l_2 c-1}{(l_2-l_1)M^{\beta_1}}  x^{\beta_1} + \frac{1- l_1 c}{(l_2-l_1)M^{\beta_2}} x^{\beta_2},$$
 and $$\ex_{x,2}[e^{-r \eta} ] = l_1 \frac{l_2 c-1}{(l_2-l_1)M^{\beta_1}}  x^{\beta_1} + l_2 \frac{1- l_1 c}{(l_2-l_1)M^{\beta_2}} x^{\beta_2}.$$
Now we choose $Z(t)= M I_{[M,\infty)\times \set{2}}(X(t),\al(t))$, $t\ge 0$. Also, let $$\wdt \eta :=\inf\set{t\ge 0: (X(t),\al(t))\in [M,\infty)\times\set{2}}.$$ Then we have
$\wdt \eta \le \eta\le \tau$.  Therefore the  fact $ 1 > \beta_2>0$ leads to
\bed \barray
J(x,1,Z) \ad = \ex_{x,1} \int_0^\tau e^{-rs} dZ(s)   \ge M \ex_{x,1} [e^{-r\wdt\eta}] \ge M \ex_{x,1} [e^{-r \eta}] \\
 \ad =   \frac{l_2 c-1}{(l_2-l_1)}  x^{\beta_1} M^{1-\beta_1} + \frac{1- l_1 c}{(l_2-l_1)}x^{\beta_2} M^{1- \beta_2}  \\
\ad \to \infty, \hbox{ as } M \to\infty.
\earray \eed
Similar calculation shows that
\bed \barray J(x,2,Z)\ad \ge \frac{l_1(l_2 c-1)}{(l_2-l_1)}  x^{\beta_1} M^{1-\beta_1} + \frac{l_2(1- l_1 c)}{(l_2-l_1)}x^{\beta_2} M^{1- \beta_2}  \\
\ad \to \infty, \hbox{ as } M \to\infty.
\earray \eed
The claim that $V(x,1)=V(x,2)=\infty $ thus follows.

Case 4: $\mu_1 \ge r$ and $\mu_2>r$. As in Case 3, we can show that the characteristic equation of (\ref{exm-2regimes})  $h(x)=g_1(x)g_2(x)-\lambda_1 \lambda_2=0$ has a solution $0< \beta_2 < 1$ and hence similar arguments as in Case 3 reveals that $V(x,1)=V(x,2)=\infty$.

\end{exm}

\begin{rem}\label{rem-optimal-control}
In Theorem \ref{thm-verification} part (b), condition \eqref{cond-op1} suggests that the optimal harvesting strategy $\wdt Z$ will harvest only in the harvesting region $$\mathcal H = \rr_+ \times \M -\mathcal C = \set{(x,\al) \in S \times \M: \phi'(x,\al) =f(x,\al)},$$
which, in turn, implies that $\phi(x,\al) = \int^x  f(y,\al)dy  $ for $(x,\al) \in \mathcal H$.
Further, if we do  harvest at time $s$ (so $\wdt X(s) \not= \wdt X(s-)$), the amount of harvest $\Delta \wdt Z(s)$ must satisfy condition \eqref{cond-op2}.
In other words, if we denote $\wdt X(s-) =x$, $\al(s-)= \alpha$, and  $\Delta \wdt Z(s) = \wdt X(s-) -\wdt X(s) =\delta x> 0$, then we must have
\begin{equation}\label{eq:about-e-optimal}  -f(x,\alpha) \delta x= \phi(x-\delta x,\alpha) -\phi(x,\alpha) = - \int^x_{x-\delta x} f(y,\alpha)dy.\end{equation}
However, if $f(x,\alpha)$ is strictly decreasing with respect to $x$ for $(x,\alpha) \in \mathcal H$, \eqref{eq:about-e-optimal} can never be satisfied. In other words, there is no {\em admissible optimal} harvesting strategy at all.
 Then a natural question arises:  Can we find an {\em admissible $\e$-optimal} harvesting policy? In the following example, the answer  to this question is positive.
\end{rem}

\begin{exm}\label{example-non-constant-price}
 As in Example \ref{example-2regimes}, let the harvested process be given by \eqref{exam-harv} and   the random environments be modeled by a two-state continuous time Markov chain $\al$ whose generator is $Q$.  Our objective is to maximize  the expected total discounted income from harvest
\begin{equation}\label{exm-non-constant-price}
V(x,\al)= \max_{Z\in \cal A} \ex_{x,\al} \int_0^\tau e^{-rs}(1+ \hat X(s-))^{-\gamma} dZ(s),
\end{equation}
where $0<\gamma<1$,  and $r$, $\tau, $ and $Z$ are as in Example \ref{example-2regimes}.

Assume that $\mu_1 > r$ and $\mu_2 > r$. As a result, the  positive roots $$p_i=\frac{1}{2}- \frac{\mu_i}{\sigma_i^2}+ \sqrt{\(\frac{1}{2}- \frac{\mu_i}{\sigma_i^2}\)^2 + \frac{2r}{\sigma_i^2}}$$ of the equations $\frac{1}{2}\sigma_i^2 x(x-1) + \mu_i x -r =0$  satisfy $0< p_i <1$, $ i=1,2$.
Suppose that \begin{equation}\label{exm2-condition} p_1=p_2 =p. \end{equation} Note that there are many nontrivial examples (in the sense that $\mu_1 \not=\mu_2$ and $\sigma_1 \not=\sigma_2$) where   condition \eqref{exm2-condition} is satisfied. For example, if $\mu_1 =1, \sigma_1^2 =2, \mu_2 = 2-\sqrt r,$ and $ \sigma_2^2 =4$, where $0< r < 1$, then
$p_1 =p_2 =p = \sqrt r$. Under condition \eqref{exm2-condition}, we   compute
\bed h(p) = g_1(p) g_2(p) -\lambda_1 \lambda_2 = (-\lambda_1) (-\lambda_2) -\lambda_1 \lambda_2 =0, \eed
where $h, g_1$, and $g_2$ are defined in Example \ref{example-2regimes}.
Consequently, it follows that $(\op-r) x^p =0$.

Assume $1-p < \gamma < 1$. Detailed calculations reveal that
\begin{equation}\label{exm2-value}
\phi(x,1)=\phi(x,2)= \begin{cases}\dfrac{(1+b)^{-\gamma}}{p b^{p-1}}x^{p}, & \text{ if } x \in (0,b),\\ \dfrac{(1+x)^{1-\gamma}-(1+b)^{1-\gamma}}{1-\gamma}+ \dfrac{ b (1+b)^{-\gamma}}{p}, & \text{ if } x \in [b, + \infty),\end{cases}
\end{equation}
 solves the quasi-variational inequality
 \bed \max\set{(\op-r) \phi(x,\al), (1+x)^{-\gamma}- \phi'(x,\al)}=0, \ \  x> 0, \ \al=1,2, \eed where $b=\frac{1-p}{p+\gamma -1}>0$. Therefore by virtue of \thmref{thm-verification}, $V(x,\al) \le \phi(x,\al)$ for all $x>0$ and $\al=1,2$.

 Next, we construct an $\e$-optimal harvesting strategy $Z^\e$.
 To this end, we denote  the  continuation region  by ${\mathcal C}:= (0,b)\times \set{1,2}  $
and the harvesting region by ${\mathcal H}:= [b,\infty) \times \set{1,2}  $.
Let $(x,\al)\in \mathcal C$ and $L_b $ be the local time process for the process
$\hat X $ at the point $b$. Then for the function $(x,\al)\mapsto x^p$, we have
\bed \ex_{x,\al} [e^{-r T_N} \hat X(T_N)^p] -x^p = \ex_{x,\al} \int_0^{T_N} e^{-rs} (\op-r) \hat X(s)^p ds
- \ex_{x,\al} \int_0^{T_N} e^{-rs} p \hat X(s)^{p-1} dL_b(s), \eed
where $T_N = \tau \wedge N \wedge \beta_N$ as in the proof of \thmref{thm-verification}.
Thus it follows that
\bed p b^{p-1} \ex_{x,\al}  \int_0^{T_N} e^{-rs} dL_b(s)= x^p - \ex_{x,\al} [e^{-r T_N} \hat X(T_N)^p].\eed
Using \eqref{beta-N}, we can readily verify that $\ex_{x,\al} [e^{-r T_N} \hat X(T_N)^p] \to 0$ as $N\to \infty$. Hence  by letting $N\to \infty$, we obtain
\begin{equation}\label{eq:local-time-mean} \ex_{x,\al}  \int_0^{\tau} e^{-rs} dL_b(s) = \frac{x^p}{p b^{p-1}}.\end{equation}
Now let $Z  =L_b$. Then \bed \barray J(x,\al, Z) \ad  = \ex_{x,\al} \int_0^\tau e^{-rs} (1+ \hat X(s))^{-\gamma} dZ(s)  = \ex_{x,\al} \int_0^\tau e^{-rs} (1+ \hat X(s))^{-\gamma} dL_b(s) \\ \ad  = (1+b)^{-\gamma} \ex_{x,\al}  \int_0^{\tau} e^{-rs} dL_b(s) = \frac{(1+b)^{-\gamma}}{p b^{p-1}}x^p = \phi(x,\al).\earray\eed

Next, let $(x,\al) \in \mathcal H$. If $x=b$, define $Z^1 \in \mathcal A$ such that  $ Z^1(0)= \Delta Z^1(0) =\varrho>0$ and  $Z^1(t)= L_b(t)$ for $t>0$.
 Then \begin{equation}\label{eq:J-Z1}\barray J(x,\al,Z^1)\ad  = (1+b)^{-\gamma} \varrho + J(b-\varrho,\al, L_b) \\ \ad =(1+b)^{-\gamma} \varrho +  \frac{(1+b)^{-\gamma}}{p b^{p-1}}(b-\varrho)^p  \\ \ad > \frac{b(1+b)^{-\gamma}}{p } -\e/3 = \phi(x,\al) -\e/3 \earray\end{equation} for $\varrho $ sufficiently small.

If $x> b$. Then as in the proof of \thmref{cor-1d-optimal},
for $n$ sufficiently large, let $\delta = \frac{x-b}{n}$, $\varsigma=n^{-5}$, and $t_i = i \varsigma/n$ for  $i=0,1,\dots, n-1$. Let $Z^2 \in \mathcal A$ increase only at times $t_i$, $i=0, 1, \dots, n-1$. More specifically, $Z^2(t_0)=\Delta Z^2(t_0) =\delta$ and   for $i=1, 2, \dots, n-1$,
$$\Delta Z(t_i) = \( \hat X(t_i-)- (x- (i+1)\delta)\)^+.$$ Denote $\theta:=\inf\set{t\ge 0: \hat X(t) \le b}$. Note that $\theta \le t_{n-1} < n^{-5}$ and  $\hat X(\theta)= b$ a.s. Define
\begin{equation}\label{eq:Z-e} Z^\e (t) = Z^2(t) I_{\set{t< \theta}} + Z^1 (t) I_{\set{t\ge \theta}}.\end{equation}

It is easy to verify that for any $x,y>0$, we have % the function $f$ is Lipschitz continuous.
\bed \abs{(1+x)^{-\gamma}-(1+y)^{-\gamma}} \le \abs{x-y},\eed
and hence $f$ is Lipschitz continuous.
As a result, similar arguments as those in the proof of \thmref{cor-1d-optimal} reveal that
\begin{equation}\label{eq:J-Z2}\barray  \ex_{x,\al}\disp\int_0^\theta e^{-r s} (1+\hat X(s-))^{-\gamma} dZ^2(s)  \ad > \int_b^x (1+y)^{-\gamma} dy -\frac{\e}{2}
    %\\[1.5ex] \ad =  \dfrac{(1+x)^{1-\gamma}-(1+b)^{1-\gamma}}{1-\gamma} - \frac{\e}{2}
    \earray \end{equation} by choosing $n$ sufficiently large.
Then for   $n$   sufficiently large,
we have from \eqref{eq:J-Z1} and \eqref{eq:J-Z2} that $$\barray J(x,\al,Z^\e)\ad  = \ex_{x,\al}\int_0^\tau e^{-r s} (1+\hat X(s-))^{-\gamma} dZ^\e(s) \\
 \ad =   \ex_{x,\al}\int_0^\theta e^{-r s} (1+\hat X(s-))^{-\gamma} dZ^2(s) +  \ex_{x,\al}\int_\theta^\tau e^{-r s} (1+\hat X(s-))^{-\gamma} dZ^1(s) \\
 \ad \ge  \ex_{x,\al}\int_0^\theta e^{-r s} (1+\hat X(s-))^{-\gamma} dZ^2(s)  + \ex_{x,\al}[ e^{-r \theta} (J(b,1, Z^1) \wedge J(b,2, Z^1))  ]  \\
 \ad > \int^x_b (1+y)^{-\gamma}dy - \e/2 + e^{-rn^4}\(\frac{b(1+b)^{-\gamma}}{p}- \e/3\) \\
 \ad >  \dfrac{(1+x)^{1-\gamma}- (1+b)^{1-\gamma}}{1-\gamma} + \frac{b(1+b)^{-\gamma}}{p} -\e =\phi(x,\al) -\e,\earray $$
This shows that the harvesting policy $Z^\e$ is $\e$-optimal.

\end{exm}

\section{Properties of the Value Function}\label{sect-properties-of-V}
  \thmref{thm-verification} gives sufficient conditions for a function $\phi$ to coincide with the value function.
  In particular, the function $\phi$ must satisfy the system of quasi-variational inequalities (\ref{quasi-vari}).
  It is natural to ask whether the converse is true: ``Does the value function $V$ defined in (\ref{value}) always satisfy (\ref{quasi-vari})?''
  In general, the answer is no since the solution $V$ is not necessarily smooth enough.
   An alternative definition for a solution to the quasi-variational inequalities (\ref{quasi-vari}) is that of a {\em viscosity solution} (see Fleming and Soner \cite{FlemingS}).

Therefore this section is devoted to the properties of the value function $V$. We present   sufficient conditions under which the value function $V$ is continuous. Also, we  show that $V$ is a viscosity solution to the coupled system of quasi variational inequalities \eqref{quasi-vari}.
\subsection{Continuity}
As we indicated in Section \ref{sect-intro}, in the definition of $J$ in \eqref{income-J-defn}, both the extinction time $\tau$ and the harvesting strategy $Z$ may depend on the initial condition $\hat X(0-)=x$.
Consequently the standard arguments for continuity of $V$ using Lipschitz continuity and Gronwall's inequality (as in \cite[Propsition 4.3.1]{YongZ}) or Dini's theorem (as in \cite[Lemma V.2.1]{FlemingS}) do not apply here. In the sequel, we establish the continuity of $V$ by first studying some elementary properties of $V$.

The following lemma is immediate, which asserts that the value function $V(x,\al)$ is nondecreasing with respect to the $x$ variable.
\begin{lem}\label{lem-increasing}
For each $\al \in \M$ and any $0<y\le x$, we have
\begin{equation}\label{value-increasing}V(x,\al)\ge   f(x,\al) ( x-y) + V(y,\al).\end{equation}
% $V(y,\al) \le V(x,\al)$.
\end{lem}
\begin{proof} Fix some $\al\in \M$.
If $y \le x$, then for any harvesting strategy $Z\in \cal A$, we define $\wdt Z$
to be the harvesting strategy such that % $\wdt Z(0) =Z(0) + (y-y)$ and that
 $\wdt Z(t)= Z(t)+ (x-y)$ for any $t\ge  0$. It is obvious that $\wdt Z \in \cal A$.
 Note also that $\Delta \wdt Z(t)= \Delta Z(t)$ for all $t > 0$.
Let $\hat X$ denote the process satisfying (\ref{harvest}) with initial condition $\hat X(0-)=y$ and harvesting strategy $Z$. Similarly, $\wdt X$ denotes the process satisfying
(\ref{harvest}) with initial condition $\wdt  X(0-)=x$ and harvesting strategy $\wdt Z$. Then we have $\hat X(t)= \wdt X(t)$ for all $t>0$. Consequently, it follows that
 $$J(x,\al, \wdt Z)= f(x,\al) ( x-y)  + J(y,\al,Z).$$
 Since $V(x,\al) \ge J(x,\al,\wdt Z)$, we  have
 $$ V(x,\al) \ge  f(x,\al) ( x-y)  + J(y,\al,Z),$$
   from which \eqref{value-increasing} follows by taking supremum over $Z\in \cal A$.
\end{proof}

\begin{lem}\label{lem-value-increase-not-much} For each $\al \in \M$
  and any $0< y \le x$, we have
 \begin{equation}\label{value-increase-not-much}V(x,\al) \le V(y,\al) +
  \max\set{V(x-y,j): j=1,\dots, m}.
\end{equation}
\end{lem}

\begin{proof}
Let $ x>y> 0$. We consider a harvested process $\hat X$   with initial conditions
$\hat X(0-) =x, \al(0) =\al$ and harvesting strategy $Z\in \cal A$. Define
$\theta  :=\inf\{t\ge 0: \hat X  (t) \le x-y\}$ and $ \tau  :=\inf\{t\ge 0: \hat X (t) \not\in S\}$.
Then $\tau \ge \theta$.

Case 1: $\theta =\infty$. Note that $\tau =\infty$.
Let $y\in (0, x)$ and  $\wdt X$ be an  another harvested process  with initial conditions
 $\wdt X(0-)=y, \al(0)=\al$ and harvesting strategy $\wdt Z \in \cal A$,
where we  choose $\wdt Z(t) = Z (t)$ for all $t\ge 0$. Consequently,
it follows that $0\le \wdt X(t) = \hat X(t) -(x-y) \le \hat X(t)$ for all $t\ge 0$.
Using the assumption that $f(\cdot, \al)$ is nonincreasing for each $\al \in \M$, we have
\beq{1d-cont-case1}\barray J(x,\al,Z) \ad = \ex\int_0^\tau e^{-rs} f(\hat X(s-),\al(s-))dZ(s) \\[1.5ex]
\ad \le  \ex\int_0^\tau e^{-rs} f(\wdt X(s-),\al(s-))d\wdt Z(s) \\
 \ad = J(y,\al,\wdt Z)\le V(y,\al).\earray \eeq

Case 2: $\theta < \infty$. Then we can write \bed\barray J(x,\al,   Z) \ad = \ex \int_0^{  \tau} e^{-rs} f(\hat X(s-),\al(s-))d  Z(s) \\
 \ad = \ex \bigg[\int_{[0,\theta)}e^{-rs} f(\hat X(s-),\al(s-)) d Z(s) + e^{-r\theta} f(\hat X(\theta-),\al(\theta-))\Delta Z(\theta)
       \\ \aad \    +
 \int_{(\theta, \tau]}e^{-rs} f(\hat X(s-),\al(s-))d   Z(s) \bigg].
 \earray\eed As in Case 1, we consider  $y\in (0, x)$ and     another harvested process $\wdt X$ with initial conditions
 $\wdt X(0-)=y, \al(0)=\al$ and harvesting strategy $\wdt Z \in \cal A$,
where we  choose $\wdt Z(t) = Z (t)$ for all
 $0\le t < \theta $   and $\wdt Z(\theta) = Z(\theta-) + \wdt X(\theta-)$. As a result,
  $\hat X(t)-\wdt X(t) =x-y$ for all $0\le t < \theta$ and  $\wdt X(\theta) = 0$.
 Then it follows from the monotonicity of $f$ that
\bed\barray J(x,\al,   Z)
  %\ad = \ex \int_0^{  \tau} e^{-rs} f(\hat X(s-),\al(s-))d  Z(s) \\
 %\ad = \ex \bigg[\int_{[0,\theta)}e^{-rs} f(\hat X(s-),\al(s-)) d Z(s)
 %+ e^{-r\theta} f(\hat X(\theta-),\al(\theta-))\Delta Z(\theta)
       %-Z(\theta-)]
 %\\ \aad \    +
 %\int_{(\theta, \tau]}e^{-rs} f(\hat X(s-),\al(s-))d   Z(s) \bigg] \\
 \ad \le \ex  \bigg[\int_{[0, \theta)}e^{-rs} f(\wdt X(s-),\al(s-)) d \wdt Z(s)
 + e^{-r\theta} f(\wdt X(\theta-),\al(\theta-))
 \Delta \wdt Z(\theta) %-\wdt Z(\theta-)]
 \\ \aad \     +\    e^{-r\theta} f(\hat X(\theta-),\al(\theta-))\!\big[ \Delta Z(\theta) -\Delta \wdt Z(\theta)
     %-(\wdt Z(\theta)-\wdt Z(\theta-))
 \big]
 %\\ \aad  \ \quad
 + \! \int_{(\theta, \tau]}\! e^{-rs} f(\hat X(s-),\al(s-))d   Z(s) \bigg ]\\
 \ad \le \ex \bigg[\int_{[0, \theta]}e^{-rs} f(\wdt X(s-),\al(s-)) d \wdt Z(s) +  \ e^{-r\theta} f(\hat X(\theta-),\al(\theta-)) \Delta Z(\theta)
         %-Z(\theta-))]
  \\ \aad \    +  \int_{(\theta, \tau]}e^{-rs} f(\hat X(s-),\al(s-))d   Z(s) \bigg]
  \\ \ad = J(y,\al, \wdt Z) + \ex  \int_{[\theta, \tau]} e^{-rs} f(\hat X(s-),\al(s-))d   Z(s).
 \earray \eed
Note that $\hat X(\theta) \le x-y$. Therefore by virtue of \lemref{lem-increasing}, we have
 \bed\barray
 \ex  \disp\int_{[\theta, \tau]} e^{-rs} f(\hat X(s-),\al(s-))d   Z(s)\ad
  = \ex\bigg[\ex\Big[\int_{[\theta, \tau]} e^{-rs} f(\hat X(s-),\al(s-))d   Z(s) \big| \cal F_\theta \Big] \bigg]
 \\[1.9ex]
 %\ad =\ex\left[ e^{-r\theta}\ex_{\hat X(\theta),\al(\theta)}\int_{[0, \tau ]} e^{-rs} f(\hat X(s-),\al(s-))d   Z(s)\right]\\
 \ad \le  \ex[e^{-r\theta} V(\hat X(\theta),\al(\theta))]\\
 \ad \le \max\set{V(x-y,j), j=1,\dots,m}.
 \earray \eed
 Hence it follows that
 \beq{1d-cont-case2}\barray J(x,\al,Z) \ad \le J(y,\al,\wdt Z) + \max\set{V(x-y,j), j=1,\dots,m} \\ \ad \le V(y,\al) + \max\set{V(x-y,j), j=1,\dots,m}.\earray\eeq

 Combining cases 1 and 2, we conclude that for any $Z\in \cal A$,
  we have
 \bed J(x,\al,Z) \le V(y,\al) + \max\set{V(x-y,j), j=1,\dots,m}.\eed
Now \eqref{value-increase-not-much} follows by taking supremum over $Z\in \cal A$.
 \end{proof}

For $h\ge 0$, we denote
\begin{equation}\label{zeta-h}
\zeta_h:= \inf\set{t\ge 0: X^{x,\al}(t)=h}, \text{ where }  X^{x,\al}(t)  \text{ is the solution to }(\ref{dyna-compo}).
\end{equation}

\begin{lem}\label{lem-V-at-0}
 Suppose that for any $t> 0$ and $h> 0$, we have
\begin{equation}\label{assump-A} \pr_{x,\al}\set{\max_{0\le s \le
      \zeta_0 \wedge t} X(s) < h} \to 1, \text{ as } x \downarrow 0,
  \text{ for each } \al\in \M.\end{equation}
If $V(x_0,\al)< \infty$ for some $x_0 >0$ and every  $\al\in  \M$, then
\begin{equation}\label{eq:V-to-0} \lim_{x\downarrow 0} V(x,\al)=0, \text{ for each } \al \in \M.\end{equation}
\end{lem}

\begin{proof} It follows from Lemmas \ref{lem-increasing} and \ref{lem-value-increase-not-much} that
$V(x,\al) < \infty$ for all $(x,\al) \in S\times \M$.
Fix some $(x,\al)\in \rr^+ \times \M$ and $Z\in \cal A$.
For any $h>0$ and  $t>0$, let $\eta:=\tau\wedge \zeta_h \wedge t$, where $\tau$ is the extinction time   defined in (\ref{tau}) and $\zeta_h$   in (\ref{zeta-h}). Then $\eta \le \tau$ and $\eta < \tau $ if and only if $\tau > \zeta_h \wedge t$.
Note that (\ref{assump-A}) implies that
$\pr_{x,\al} \set{\zeta_0 < \zeta_h \wedge t} \to 1,$ as $x \downarrow 0$, for each $\al \in \M$. Hence for any $\e >0$, we can find a $\delta >0$ such that
\beq{prob=1} \pr_{x,\al}\set{\zeta_0 < \zeta_h \wedge t} > 1-\e, \text{ for any }0\le x< \delta \text{ and each }\al\in \M.\eeq
Note that for any $s\ge 0$, $\hat X^{x,\al}(s)\le X^{x,\al}(s)$ by the admissibility of $Z\in \cal A$. Hence it follows that $\tau \le \zeta_0$.
This, together with  \eqref{prob=1}, implies  that $\pr_{x,\al} \set{\tau <\zeta_h \wedge t} > 1-\e.$ That is
\beq{prob=0} \pr_{x,\al} \set{\tau > \eta} < \e,\text{ for any }0\le x< \delta \text{ and each }\al\in \M.\eeq
On the other hand, by virtue of  \cite[Proposition 2.3]{YZ-10}, we have
\beq{E=0} \ex_{x,\al} \left[\sup_{0 \le s \le \zeta_0 \wedge t} X(s)\right] < \e, \text{ for any }0\le x< \delta \text{ and each }\al\in \M.\eeq
Now we compute
\bed \begin{aligned}
J(x,\al,Z)& = \ex_{x,\al} \int_0^\tau e^{-r s} f(\hat X(s-),\al(s-))dZ(s) \\
& = \ex_{x,\al} \left[\int_0^\eta  e^{-r s} f(\hat X(s-),\al(s-))dZ(s)  \right] \\ & \quad + \ex_{x,\al} \left[I_{\set{\tau> \eta}}\int_\eta^\tau  e^{-r s} f(\hat X(s-),\al(s-))dZ(s)  \right]\\ &: = A + B.\end{aligned}\eed
Since $f(\cdot,\al)$ is nonincreasing and that $\hat X(s-) \ge 0$ for all $0\le s \le \eta $, it follows that $$f(\hat X(s-),\al(s-)) \le f(0, \al(s-)) \le \max_{i\in \M} f(0,i) =K,$$ where in the above and hereafter, $ K$ is a generic positive constant not depending
on $x$ or $t$ whose exact value may change in different appearances. Hence it follows that
$ A \le K \ex_{x,\al}[Z(\eta)]$.
Further, since $Z$ is an admissible harvesting strategy, we have
$$Z(\eta) \le X(\eta) \le \sup_{0\le s \le \zeta_0 \wedge t} X(s).$$
 Note in the above, we  used the fact that $\eta = \tau\wedge \zeta_h \wedge t \le \tau \wedge t \le \zeta_0 \wedge t$.
Thus we have from (\ref{E=0}) that $A \le \ex_{x,\al}[Z(\eta)] < K \e.$

On the other hand,
 \bed \begin{aligned} B & = \ex_{x,\al} \left[I_{\set{\tau> \eta}}\int_\eta^\tau  e^{-r s} f(\hat X(s-),\al(s-))dZ(s)  \right]
 \\ & = \ex_{x,\al} \left[I_{\set{\tau> \eta}} \ex\Big[\int_\eta^\tau  e^{-r s} f(\hat X(s-),\al(s-))dZ(s)| \F_{\eta}\Big]\right]
  \\ & \le \ex_{x,\al}\left[I_{\set{\tau> \eta}} e^{-e \eta}V(X(\eta),\al(\eta)) \right].  \end{aligned}\eed
By the definition of $\eta$, we have $\hat X(\eta) \le X(\eta) \le h$.   Thus \lemref{lem-increasing} leads to
$$V(X(\eta),\al(\eta)) \le V(h,\al(\eta)) \le \max_{i\in \M}V(h,i) < \infty.$$
Therefore it follows from the above computations and (\ref{prob=0}) that
$$J(x,\al,Z) \le K\e + \max_{i\in\M} V(h,i) \pr_{x,\al}\set{\tau> \eta}  < K \e.$$
Taking supremum over $Z\in \cal A$, we obtain $V(x,\al) < K\e$, for all $0\le x< \delta$  and each $\al\in \M$. Therefore \eqref{eq:V-to-0} follows and this completes the proof of the lemma.
\end{proof}

By virtue of Lemmas \ref{lem-increasing}, \ref{lem-value-increase-not-much}, and \ref{lem-V-at-0}, we have the following theorem, which presents a sufficient condition for continuity of the value function.
\begin{thm}\label{thm-value-cont}
Let the conditions of \lemref{lem-V-at-0} be satisfied.
Then the value function $V$ defined in \eqref{value} is continuous with respect to the variable $x$.
\end{thm}

\begin{rem}\label{about-assump-A}
Note that (\ref{assump-A}) is the crucial assumption in \thmref{thm-value-cont}, it also plays a key role of the proof of \lemref{lem-V-at-0}. One may wonder under what condition(s), is (\ref{assump-A}) valid? \begin{exm}
If the unharvested process is given by
 \beq{1d-special} dX(t)= b(\al(t)) dt+ \sigma(\al(t)) dw(t), \eeq
 where for each $\al \in \M$, $b(\al)\in \rr$ and $\sigma(\al) >0$, then
 using the same argument as that of \cite[Lemma 1]{Choulli}, we obtain \eqref{assump-A}.
\end{exm}
Next we present a sufficient condition for (\ref{assump-A}). \end{rem}

\begin{prop} \label{l-W}If there exists a function $W:S\times \M
  \mapsto \rr^+$ satisfying
\begin{itemize} \item[{\em (i)}] for each $\al\in\M$, $W(\cdot,\al)$ is continuous on $[0,\infty)$ and vanishes only at $x=0$,
\item[{\em (ii)}] $\op W(x,\al) \le 0$ for all $(x,\al)\in S\times \M$.
\end{itemize}
Then for any $t> 0$ and $h> 0$, we have
\bed \pr_{x,\al}\set{\max_{0\le s \le \zeta_0 \wedge t} X(s) < h} \to 1, \text{ as } x \downarrow 0, \text{ for each } \al\in \M.\eed
\end{prop}
\begin{proof}
Fix some $t>0$ and $h>0$. Let $(x,\al) \in S \times \M$ with $x<h$. Denote $W_h:=\inf\set{W(y,j): y\in S, y\ge h, j\in\M}$. Then assumption (i) implies that
$W_h >0$. By virtue of It\^o's formula and assumption (ii), we have
\bed \barray\ex_{x,\al} W(X(t\wedge \zeta_0\wedge \zeta_h), \al(t\wedge \zeta_0\wedge \zeta_h)) \ad = W(x,\al) + \ex_{x,\al} \int_0^{t\wedge \zeta_0\wedge \zeta_h}
\op W(X(s),\al(s))ds  \\
\ad\le W(x,\al).\earray\eed
Now since $W$ is nonnegative, it follows that
\bed W(x,\al) \ge \ex_{x,\al} \left[W(X(\zeta_h),\al(\zeta_h))I_{\set{\zeta_h < \zeta_0\wedge t}}\right] \ge W_h \pr_{x,\al}\set{\zeta_h < \zeta_0\wedge t}.\eed
Thus we have
\bed \pr_{x,\al}\set{\zeta_h < \zeta_0\wedge t} \le \frac{W(x,\al)}{W_h}.\eed
This, together with assumption (i), leads to
 \bed \pr_{x,\al}\set{\zeta_0\wedge t < \zeta_h} \ge 1- \frac{W(x,\al)}{W_h}  \to 1, \ \ \hbox{ as }x \downarrow 0.\eed
Note that $\zeta_0\wedge t < \zeta_h$ if and only if $\max_{0\le s \le \zeta_0\wedge t}X(s) < h$. Therefore the desired assertion follows.
\end{proof}

\subsection{Viscosity Solution}
In this subsection, we aim to characterize the value function as a
viscosity solution of the   coupled system of quasi-variational inequalities \eqref{quasi-vari}.
Let's  recall the notion of viscosity solution.

\begin{defn}
  \label{d-visc}  {\rm A function
$u$ is said to be
a {\em viscosity subsolution}   of
\eqref{quasi-vari}, if
 for any $(x_0,\al_0) \in S \times \M$ and function
  $\varphi(\cdot, \al) \in C^2(S)$ satisfying $\varphi(\xalz) = u(\xalz)$ and $\varphi(x,\al) \le u(x,\al)$ for
  all $x$ in a neighborhood of $x_0$ and each $\al \in \M$, we have
  $$ \max\set{(\op-r) \varphi(x_0,\al_0),  f(\xalz)-  \varphi'(\xalz)} \le 0. $$
 Similarly, a function
$u$ is said to be
a {\em viscosity supersolution}   of
\eqref{quasi-vari}, if
 for any $(x_0,\al_0) \in S \times \M$ and function
  $\varphi(\cdot, \al) \in C^2(S)$ satisfying $\varphi(\xalz) = u(\xalz)$ and $\varphi(x,\al) \ge u(x,\al)$ for
  all $x$ in a neighborhood of $x_0$ and each $\al \in \M$, we have
  $$ \max\set{(\op-r) \varphi(x_0,\al_0),  f (\xalz)- \varphi'(\xalz)} \ge 0. $$
The function
$u$ is said to be a {\em viscosity solution}  of
\eqref{quasi-vari}, if it
is both a viscosity
subsolution and a viscosity supersolution.
}\end{defn}

\begin{thm}\label{thm-viscosity}
Assume the conditions of \thmref{thm-value-cont}.
Then  the value function
$V$ is a viscosity solution of the coupled system of quasi-variational inequalities \eqref{quasi-vari}
with boundary condition \eqref{eq:V-to-0}.
 \end{thm}

 \begin{proof} The proof is motivated by
   \cite{Choulli}  and
   \cite[Theorem VIII 5.1]{FlemingS}.
 We divide the proof into two parts. The first part shows that $V$ is a viscosity subsolution of (\ref{quasi-vari}), while the second part establishes that $V$ is viscosity supersolution of (\ref{quasi-vari}).

 Step 1. We show that $V$ is a viscosity subsolution of (\ref{quasi-vari}). That is, for any $(x_0,\al_0) \in S\times \M$ and any $C^2$ function $\phi(\cdot, \cdot)$ satisfying $\phi(x_0,\al_0)= V(x_0,\al_0)$ and that $\phi(x,\al) \le V(x,\al)$ for all $x$ in a neighborhood of $x_0$ and each $\al\in \M$, we have
 \begin{equation}\label{viscosity-sub}
 \max\set{(\op-r) \phi(\xalz), f(\xalz)-  \phi'(x_0,\al_0)} \le 0.\end{equation}
 Let $B_\e(x_0):=\set{x\in \rr: \abs{x-x_0} < \e}$, where $\e>0$ is sufficiently small so that (i) $\lbar B_\e(x_0) \subset S$ and (ii) $\phi(x,\al)
\le V(x,\al)$ for all $(x,\al)\in \lbar B_\e(x_0)\times \M$, where $\lbar B_\e(x_0)= \set{x\in \rr: \abs{x-x_0}\le \e}$ denotes the closure of $ B_\e(x_0)$.
Choose $Z\in \cal A$ such that $Z(0-)=0$ and $Z(t)=\eta $ for all $t \ge 0$, where $0\le \eta < \e$.
Let $\hat X\cd= \hat X^{\xalz}\cd$ be the corresponding harvested process with initial condition $(\xalz)$ and harvesting strategy $Z\cd$.
Put \bed \theta:= \inf\set{t \ge 0: \hat X(t) \notin B_{\e }(x_0)}.\eed
Note that the chosen harvesting strategy $Z$ guarantees that $\hat X\cd$ has at most one jump at $t=0$ and remains continuous on $(0,\theta]$.
This, together with  the choice of $\e$, implies that  $\theta \le \tau$ and that
$\hat X(t) \in \lbar B_\e(x_0)$ for all $0\le t \le \theta$. By virtue of the dynamic programming principle (\ref{dyna-prog}), for any $h>0$, we have
\beq{sub-ineq1} \begin{aligned} \phi(\xalz) & =V(\xalz) \\ & \ge \ex \bigg[ \int_0^{\theta\wedge h} e^{-rs}   f(\hat X(s-),\al(s-))  dZ(s)  + e^{-r ( \theta\wedge h)} V(\hat X(\theta\wedge h),\al(\theta\wedge h))\bigg] \\
& \ge \ex \bigg[ \int_0^{\theta\wedge h} e^{-rs}  f(\hat X(s-),\al(s-)) dZ(s) +  e^{-r (\theta\wedge h)} \phi(\hat X(\theta\wedge h),\al(\theta\wedge h))\bigg]. \end{aligned}\eeq
Applying the generalized It\^o formula to the process $e^{-rs} \phi(\hat X(s),\al(s))$, we obtain
\bed \begin{aligned} & e^{-r(\theta\wedge h)} \phi(\hat X(\theta\wedge h),\al(\theta\wedge h)) - \phi(\xalz)
\\&\  = \int_0^{\theta\wedge h} e^{-r s} (\op-r) \phi(\hat X(s),\al(s))ds
+ \int_0^{\theta\wedge h} e^{-rs}    \phi'(\hat X(s),\al(s))  \sigma(\hat X(s),\al(s))dw(s)  \\
& \ \   - \int_0^{\theta\wedge h} e^{-rs}   \phi'(\hat X(s),\al(s)) dZ^c(s)  + \sum_{0\le s \le \theta\wedge h} e^{-rs}\left[\phi(\hat X(s),\al(s-))- \phi(\hat X(s-),\al(s-))\right].\end{aligned}\eed
Since $\phi \in C^2$, $\sigma$ is continuous, and $\hat X(s) \in \lbar B_\e(x_0)$ for all $0\le s \le \theta$, it follows that
$$ \ex\int_0^{\theta\wedge h} e^{-rs}  \phi'(\hat X(s),\al(s)) \sigma(\hat X(s),\al(s))dw(s) =0.$$
Consequently, we have
\beq{sub-eq2}\begin{aligned}\phi(\xalz)
= &\ \ex  e^{-r(\theta\wedge h)} \phi(\hat X(\theta\wedge h),\al(\theta\wedge h))
  -\ex\int_0^{\theta\wedge h} e^{-r s} (\op-r) \phi(\hat X(s),\al(s))ds \\
 & \ + \ex  \int_0^{\theta\wedge h} e^{-rs}  \phi'(\hat X(s),\al(s)) dZ^c(s)
 \\& \  - \ex  \sum_{0\le s \le \theta\wedge h} e^{-rs}\left[\phi(\hat X(s),\al(s-))- \phi(\hat X(s-),\al(s-))\right].  \end{aligned}\eeq
 A combination of (\ref{sub-ineq1}) and (\ref{sub-eq2}) leads to
 \beq{sub-ineq3}\begin{aligned}
0 \ge & \ \ex  \int_0^{\theta\wedge h} e^{-rs}  f(\hat X(s-),\al(s-))  dZ(s)
  +  \ex \int_0^{\theta\wedge h} e^{-r s} (\op-r) \phi(\hat X(s),\al(s))ds
 \\ & \ -\ex  \int_0^{\theta\wedge h} e^{-rs}  \phi'(\hat X(s),\al(s)) dZ^c(s) \\ & \ +\ex  \sum_{0\le s \le \theta\wedge h} e^{-rs}\left[\phi(\hat X(s),\al(s-))- \phi(\hat X(s-),\al(s-))\right].
  \end{aligned}\eeq
  Now let $\eta=0$, i.e., $Z(t)\equiv 0$ for any $t\ge 0$. Then (\ref{sub-ineq3}) can be rewritten as
  \bed  \begin{aligned}0  & \ge  \ex \int_0^{\theta\wedge h} e^{-r s} (\op-r) \phi(\hat X(s),\al(s))ds
    \\ & =  \ex \int_0^{ h} e^{-r s} (\op-r) \phi(\hat X(s),\al(s))I_{\set{s \le \theta}}ds.
   \end{aligned}\eed
  Note that $ e^{-r s} (\op-r) \phi(\hat X(s),\al(s))I_{\set{s \le \theta}}$ is bounded for all $0\le s\le h$ by the definition of $\theta$.
  Hence there exists some $0\le \xi^h \le h$ such that
  \bed   0    \ge \ex \int_0^{ h} e^{-r s} (\op-r) \phi(\hat X(s),\al(s))I_{\set{s \le \theta}}ds \ge   h \ex [ e^{-r\xi^h} (\mathcal{L} - r) \phi (\hat X(\xi^h),
  \alpha(\xi^h)) I_{\{\xi^h\le \theta\}}].\eed
   Note that as $h\downarrow 0$, $\xi^h \downarrow 0$.
   This implies that   $(\hat
X(\xi^h), \alpha(\xi^h)) \to (x_0, \alpha_0)$ a.s. due to the
choice of $Z \equiv 0$. With the continuity of $(\mathcal{L} -r) \phi$,
we conclude that
$$e^{-r\xi^h} (\mathcal{L} - r) \phi (\hat
X(\xi^h),  \alpha(\xi^h)) I_{\{\xi^h\le \theta\}} \to (\mathcal{L}
-r) \phi(x_0, \alpha_0), \hbox{ a.s.  as } h\downarrow  0. $$
 Therefore it follows from the bounded convergence theorem that
  \beq{sub-ineq4} (\op-r) \phi(\xalz) \le 0.\eeq
  On the other hand, if we choose   $0< \eta < \e$,  then (\ref{sub-ineq3})
  reduces to
  \bed  \ex \int_0^{\theta\wedge h} e^{-r s} (\op-r) \phi(\hat X(s),\al(s))ds  + f (\xalz) \eta  + \phi(x_0-\eta,\al_0) -\phi(\xalz) \le 0.\eed
  Now sending $h\downarrow 0$, we have \bed f (\xalz) \eta  + \phi(x_0-\eta,\al_0) -\phi(\xalz) \le 0.\eed
  Finally, divide the above inequality by $\eta $ and let $\eta  \to 0$, it follows that
  \beq{sub-ineq5} f (\xalz) -  \phi'(\xalz) \le 0. \eeq
  Now (\ref{viscosity-sub}) follows from a combination of (\ref{sub-ineq4}) and (\ref{sub-ineq5}).

  Step 2. We need to show that $V$ is also a viscosity supersolution of (\ref{quasi-vari}).
   That is, for any $(\xalz) \in S \times \M $ and any $\varphi \in C^2$ such that $\varphi(\xalz)=V(\xalz)$ and that $\varphi(\xalz)\ge V(\xalz)$ for $x$ in a neighborhood of $x_0$ and each $\al\in \M$,
   we have
   \begin{equation}\label{viscosity-super}
   \max\set{(\op-r) \varphi(\xalz), f(\xalz)-\varphi'(\xalz)} \ge 0.
   \end{equation}
 Suppose on the contrary that (\ref{viscosity-super})   was wrong, then there would exist some $(\xalz)\in S\times \M$, a $\varphi\in C^2$, and a constant $A>0$ such that
 \begin{equation}\label{vis-super-contra}
 \max\set{(\op-r) \varphi(\xalz),  f(\xalz)- \varphi'(\xalz)} \le -2A < 0.
 \end{equation}

 In what follows, we will derive a contradiction to (\ref{vis-super-contra}). This is achieved in several steps. First we use the generalized It\^o formula and (\ref{vis-supersoln-contra}) to obtain (\ref{super-ineq1}). Next, detailed analysis using the monotonicity of the functions $V$ and $f$ leads
 to (\ref{super-ineq5}). Then we claim in (\ref{super-claim2}) that the last term in (\ref{super-ineq5}) is bounded below by a positive constant, from which, with
 the aid of dynamic programming (\ref{dyna-prog}), we obtain a contradiction to (\ref{vis-super-contra}).  The final step of the proof is devoted to the proof of (\ref{super-claim2}).

 Fix some $Z\in \cal A$ and let $\hat X\cd =\hat X^{\xalz}\cd$ be the corresponding harvested process. Define $B_\e(x_0)$ as in Step 1, where
  $\e>0$ is small enough so that (i) $\lbar B_\e(\xz) \subset S$,
  (ii) $ \varphi(x,\al)\ge V(x,\al)$ for all $(x,\al) \in \lbar B_\e(x_0)\times \M$, and (iii)
   \begin{equation}\label{vis-supersoln-contra}
 \max\set{(\op-r) \varphi(x,\al), f (x,\al)- \varphi'(x,\al)} \le - A < 0, \ \ (x,\al) \in \lbar B_\e(x_0)\times \M.
 \end{equation}

 Let $\theta:=\inf\set{t\ge 0: \hat X(t) \notin B_\e(\xz)}$. Then $\theta\le \tau$. It follows from the generalized It\^o formula that
 \bed \begin{aligned}
 & \ex e^{-r\theta}\varphi(\hat X(\theta-),\al(\theta-)) - \varphi(\xalz) \\ & \ =\ex \int_0^{\theta-} e^{-rs} (\op-r) \varphi(\hat X(s),\al(s))ds
 -\ex \int_0^{\theta-} e^{-rs}   \varphi'(\hat X(s),\al(s))dZ^c (s)\\ & \ \ \ \  + \ex\sum_{0\le s < \theta}  e^{-rs}\left[\varphi(\hat X(s),\al(s-))-\varphi(\hat X(s-),\al(s-))\right].
  \end{aligned}\eed
  Note that \bed \begin{aligned}& \varphi(\hat X(s),\al(s-))-\varphi(\hat X(s-),\al(s-)) \\ & \ \ =   (\hat X (s)-\hat X (s-)) \varphi'(\hat X(s-)+z(\hat X(s)-\hat X(s-)),\al(s-))\\ &\ \ = - \Delta  Z (s)  \varphi'(\hat X(s-)+z(\hat X(s)-\hat X(s-)),\al(s-)) \end{aligned} \eed for some $z\in[0,1]$.
  But by virtue of (\ref{vis-supersoln-contra}), for all $0\le s < \theta$, we have
  \bed - \varphi'(\hat X(s-)+z(\hat X(s)-\hat X(s-)),\al(s-)) \le - f  (\hat X(s-)+z(\hat X(s)-\hat X(s-)),\al(s-))-A.\eed
  Further, since $\hat X(s) \le \hat X(s-)+z(\hat X(s)-\hat X(s-)) \le \hat X(s-)$ and   $f (\cdot,\al)$ is non-increasing, we have
  \bed - f  (\hat X(s-)+z(\hat X(s)-\hat X(s-)),\al(s-)) \le - f (\hat X(s-),\al(s-)).\eed
  Hence it follows from (\ref{vis-supersoln-contra}) that
  \bed \begin{aligned}
 & \ex e^{-r\theta}\varphi(\hat X(\theta-),\al(\theta-)) - \varphi(\xalz)
 \\ & \ \le\ex \int_0^{\theta-}  e^{-rs} (-A)ds
 +\ex \int_0^{\theta-} e^{-rs}   (-f (\hat X(s),\al(s))-A)dZ^c (s)
 \\ & \ \ \ \  + \ex\sum_{0\le s < \theta}  e^{-rs} (-f (\hat X(s),\al(s))-A) \Delta Z (s)
 \\ & = -\ex  \int_0^{\theta-}  e^{-rs}   f(\hat X(s),\al(s))  dZ(s)
% \\ & \ \ \ \
- A \ex \int_0^{\theta}  e^{-rs} ds
 - A \ex \int_0^{\theta-} e^{-rs}   dZ(s).
   \end{aligned}\eed
  Therefore
  \beq{super-ineq1}\begin{aligned} \varphi(\xalz) \ge & \ \ex e^{-r\theta}\varphi(\hat X(\theta-),\al(\theta-))
     + \ex  \int_0^{\theta-}  e^{-rs}  f(\hat X(s-),\al(s-))  dZ(s)  \\ & + A \ex \left[\int_0^{\theta}  e^{-rs} ds
 + \int_0^{\theta-} e^{-rs}   dZ(s) \right].\end{aligned} \eeq
 Note that $\hat X(\theta) \le \hat X(\theta-)$ and $\hat X(\theta-) \in B_\e(\xz)$. Thus there exists some $\lambda\in [0,1]$ such that
  $$x_\la:= \hat X(\theta-) + \la (\hat X(\theta)-\hat X(\theta-)) = \hat X(\theta-) - \la  \Delta  Z(\theta) \in \partial B_\e(\xz).$$
  Moreover, $\hat X(\theta) \le x_\la \le \hat X(\theta-)$.
  Note that $$ \begin{aligned}\varphi(\hat X(\theta-),\al(\theta-)) - \varphi(x_\la,\al(\theta-)) & =
          (\hat X(\tha-)-x_\la) \varphi'(\hat X(\tha- )+ z(\hat X(\tha)-x_\la),\al(\tha-))  \\ & = \la  \Delta Z(\tha) \varphi'(\hat X(\tha- )+ z(\hat X(\tha)-x_\la),\al(\tha-))  .\end{aligned}$$
  But (\ref{vis-supersoln-contra}) and the monotonicity of $f (\cdot, \al)$ imply that
   $$\begin{aligned}   \varphi'(\hat X(\tha- )+ z(\hat X(\tha)-x_\la),\al(\tha-)) & \ge f(\hat X(\tha- )+ z(\hat X(\tha)-x_\la),\al(\tha-)) +A\\ & \ge f (\hat X(\tha- ),\al(\tha-))+A.\end{aligned}$$
   This, together with the fact that $ \Delta Z (\tha)\ge 0$, leads to
   \beq{super-ineq2} \varphi(\hat X(\theta-),\al(\theta-)) - \varphi(x_\la,\al(\theta-)) \ge \la   \Delta  Z(\tha) \left[ f(\hat X(\tha- ),\al(\tha-))+A \right].\eeq
   Combing (\ref{super-ineq1}) and  (\ref{super-ineq2}), we obtain
   \beq{super-ineq3}\begin{aligned}
   V(\xalz) & = \varphi(\xalz)
   \\ & \ge  \ex  \int_0^{\theta-}  e^{-rs}   f(\hat X(s-),\al(s-))  dZ(s)
     + \ex e^{-r\tha} \varphi(x_\la, \al(\tha-))
   \\ &\ \  + A \ex      \left[\int_0^{\theta}  e^{-rs} ds
 + \int_0^{\theta-} e^{-rs}   dZ(s) \right] % \\ & \ \
 +   \la \ex e^{-r\tha}    \Delta Z(\tha) \left[ f(\hat X(\tha- ),\al(\tha-))+A\right].
       \end{aligned} \eeq

     Note that $x_\la \in \lbar B_\e(\xz)$. Therefore we have $\varphi(x_\la, \al(\tha-)) \ge V(x_\la, \al(\tha-))$.
   On the other hand, since $\hat X(\tha)\le x_\la$, it follows  from (\ref{value-increasing}) that
   \beq{super-ineq4}\begin{aligned}
       V(x_\la,\al(\theta-))&  \ge V(\hat X(\theta),\al(\theta-))+  [x_\la-\hat X(\tha)] f(x_\la,\al(\theta-))
    \\ & \ge V(\hat X(\theta),\al(\theta-)) + (1-\la)   \Delta Z(\theta)  f(\hat X(\tha-),\al(\theta-)).
   \end{aligned} \eeq
Note that   \beq{super-claim1}\ex e^{-r\tha} V(\hat X(\theta), \al(\tha-)) = \ex e^{-r\tha} V(\hat X(\theta), \al(\tha)).\eeq
In fact, by virtue of  \cite[Theorem 2.12]{YZ-10}, $\al\cd$ is continuous in mean square. Hence it follows that
\bed \begin{aligned}
&\abs{\ex e^{-r\tha} V(\hat X(\theta), \al(\tha-)) - \ex e^{-r\tha} V(\hat X(\theta), \al(\tha))} \\
& \ \ \le \ex \abs{e^{-r\tha}[V(\hat X(\theta), \al(\tha)) - V(\hat X(\theta), \al(\tha-)) ]I_{\set{\al(\tha) \not=\al(\tha-)}} } \\
& \ \ \le K \pr \set{\al(\tha) \not=\al(\tha-)}=0,\end{aligned}\eed
where $K$ is some positive constant. Therefore (\ref{super-claim1}) follows.
   Now put (\ref{super-ineq4}) and (\ref{super-claim1}) into  (\ref{super-ineq3}) and  we obtain
    \beq{super-ineq5}\begin{aligned}
   V(\xalz) &  \ge  \ex  \int_0^{\theta-}  e^{-rs}  f(\hat X(s-),\al(s-))  dZ(s)
     + \ex e^{-r\tha} V(\hat X(\theta), \al(\tha))
 \\ &\ \  + A \ex      \left[\int_0^{\theta}  e^{-rs} ds  + \int_0^{\theta-} e^{-rs}   dZ(s) \right]  +  (1-\la)\ex e^{-r\tha}   \Delta Z(\theta)  f(\hat X(\tha-),\al(\theta-))
\\ & \ \ +   \la \ex e^{-r\tha}  \Delta   Z(\tha)\left[ f(\hat X(\tha- ),\al(\tha-))+A\right]
% \\ & \ \  +  (1-\la)\ex e^{-r\tha}   \Delta Z(\theta)  f(\hat X(\tha-),\al(\theta-))
  \\ & =  \ex  \int_0^{\theta}  e^{-rs}   f(\hat X(s-),\al(s-))  dZ(s)  + \ex e^{-r\tha} V(\hat X(\theta), \al(\tha))
 \\ & \ \ +   A \ex      \left[\int_0^{\theta}  e^{-rs} ds
 + \int_0^{\theta-} e^{-rs}  dZ(s)  + \la  e^{-r\theta}  \Delta Z(\tha) \right].
    \end{aligned} \eeq

  We   now claim that for  some constant $\kappa >0 $, we have
  \begin{equation}\label{super-claim2}\ex      \left[\int_0^{\theta}  e^{-rs} ds
 + \int_0^{\theta-} e^{-rs}   dZ(s)  + \la   e^{-r\theta}  \Delta Z(\tha) \right] \ge \kappa. \end{equation} .

   Assume   (\ref{super-claim2}) for the moment. Then  (\ref{super-ineq5}) can be rewritten as
   \beq{super-ineq6}\begin{aligned}
   V(\xalz) &  \ge  \ex  \int_0^{\theta}  e^{-rs}   f(\hat X(s-),\al(s-))  dZ(s)
     + \ex e^{-r\tha} V(\hat X(\theta), \al(\tha)) + A \kappa.  \end{aligned} \eeq
     Taking supremum over $Z\in \cal A$, it follows that
      \beq{super-ineq7}\begin{aligned}
   V(\xalz) &  \ge \sup_{Z\in \cal A} \ex  \left[\int_0^{\theta}  e^{-rs}  f(\hat X(s-),\al(s-))  dZ(s)
     +   e^{-r\tha} V(\hat X(\theta), \al(\tha))\right] + A \kappa.  \end{aligned} \eeq
     But in view of the dynamic programming principle (\ref{dyna-prog}), \eqref{super-ineq7} can be rewritten as
     \bed  V(\xalz) \ge  V(\xalz) + A\kappa >  V(\xalz).\eed
     This is a contradiction. So we must have (\ref{viscosity-super}) and hence $V$ is a viscosity supersolution of (\ref{quasi-vari}).

     Now it remains to show  (\ref{super-claim2}). To this end, we consider the function $\wdt W(x,\al):= \abs{x-x_0}^2 -\e^2$ for $(x,\al)\in B_\e(\xz)\times \M$. Then it follows that
     \bed (\op-r) \wdt W(x,\al)= 2 ( x-\xz) b(x,\al)+ \frac{1}{2} 2  \sigma^2(x,\al) -r ( \abs{x-x_0}^2 -\e^2).\eed
 Since $\wdt W$, $b$, and $\sigma$ are continuous, and $\M$ is a finite set, it is obvious  that
 $$\abs{ (\op-r) \wdt W(x,\al)} \le K < \infty$$ for some positive constant $K$. Now let $K_0:=\frac{1}{2\e+ K}$ and define
 $W(x,\al)= K_0\wdt W(x,\al)$ for  $(x,\al)\in B_\e(\xz)\times \M$.
 Then it follows immediately that \beq{claim2-ineq1}\abs{ (\op-r)   W(x,\al)} < 1, \ \ (x,\al)\in B_\e(\xz)\times \M. \eeq
 Moreover, we have
 \beq{claim2-ineq2}   W'(x,\al) = 2K_0 ( x-\xz)  \ge -1.\eeq
  Now apply the generalized It\^o formula to $e^{-rs}W(\hat X(s),\al(s))$,
 \beq{claim2-eq}\begin{aligned} \ex& [e^{-r\theta} W(\hat X(\tha-),\al(\tha-))] - W(\xalz) \\& = \ex\int_0^{\tha-} e^{-rs} (\op-r) W(\hat X(s),\al(s))ds - \ex\int_0^{\tha-}  e^{-rs}     W'(\hat X(s),\al(s))dZ^c(s)
\\ & \ \ \ + \ex \sum_{0\le s < \tha}  e^{-rs} [W(\hat X(s),\al(s-))-W(\hat X(s-),\al(s-))]. \end{aligned}\eeq
But by virtue of (\ref{claim2-ineq2}), we have
\beq{claim2-ineq3}\begin{aligned} W&(\hat X(s),\al(s-))-W(\hat X(s-),\al(s-))
\\ & =    W'(\hat X(s-)+ z(\hat X(s)-\hat X(s-)),\al(s))(\hat X (s)-\hat X (s-))\\
& =  -    W'(\hat X(s-)+ z(\hat X(s)-\hat X(s-)),\al(s)) \Delta  Z (s)
\\ & \le  -   \Delta Z (s).\end{aligned}\eeq
Hence it follows from (\ref{claim2-ineq1})--(\ref{claim2-ineq3}) that
\beq{claim2-ineq4} \begin{aligned}\ex& [e^{-r\theta} W(\hat X(\tha-),\al(\tha-))] - W(\xalz)
 \\ & \le \ex \int_0^{\tha} e^{-rs}ds + \ex\int_0^{\tha-}  e^{-rs}   dZ^c(s) + \ex \sum_{0\le s < \theta} e^{-rs}  \Delta   Z(s) \\
& = \ex \int_0^{\tha} e^{-rs}ds + \ex\int_0^{\tha-}   e^{-rs} d Z(s). \end{aligned}\eeq

Also, recall that $\hat X(\tha) \le x_\la \le \hat X(\tha-)$. It follows from (\ref{claim2-ineq2}) that
\beq{claim2-ineq5} \begin{aligned}
 W&(\hat X(\tha-),\al(\tha-)) - W(x_\la,\al(\tha-))\\ & =   W' (x_\la + z(\hat X(\tha-)-x_\la),\al(\tha-)) \left[ \hat X(\tha-)-x_\la\right]
\\ &=\la    W' (x_\la + z(\hat X(\tha-)-x_\la),\al(\tha-)) \Delta Z(\tha) \\ &  \ge -\la  \Delta   Z(\tha).
\end{aligned}\eeq
Combining (\ref{claim2-ineq4}) and (\ref{claim2-ineq5}), we have
\bed  \ex \int_0^{\tha} e^{-rs}ds +\ex\int_0^{\tha-}   e^{-rs} d  Z(s)  + \la \ex e^{-r\theta}  \Delta  Z(\tha)
\ge \ex e^{-r\tha}W(x_\la,\al(\tha-)) -W(\xalz). \eed
But $x_\la \in \partial B_\e(\xz)$. Consequently, $W(x_\la,\al(\tha-))=0$. Also, it is immediate that $W(\xalz)=-K_0\e^2$. Hence it follows that
\bed  \ex \int_0^{\tha } e^{-rs}ds +\ex\int_0^{\tha-}   e^{-rs}  d Z(s)   + \la \ex  e^{-r\theta}  \Delta  Z(\tha)  \ge K_0\e^2 =\kappa >0.\eed
This establishes (\ref{super-claim2}) and hence finishes the proof of the theorem.
  \end{proof}

\section{Conclusions and Remarks}\label{sect-conclusion}

In this work, we considered the optimal harvesting problem for a single species living in random environments. We first established a verification theorem, based on which we explicitly constructed an $\e$-optimal harvesting strategy under additional conditions.
Next we obtained the continuity of the value function and further characterized it as a viscosity  solution of the
coupled system of quasi-variational inequalities \eqref{quasi-vari}.

In examples \ref{example-2regimes} (cases 1 and 2) and \ref{example-non-constant-price},
 thanks to the special structures of the harvesting and continuation regions,  we
were able to obtain the value functions and ($\e$-)optimal harvesting
policies. It will be very interesting to investigate whether ($\e$-)optimal harvesting
policies exist in more general settings. % A natural question is: can we obtain a general result?

The next logical step is to consider optimal harvesting strategy for multiple but finite number of interacting species in random environments.
For virtually all ecosystems, the species  often interact with each other
and a small change of one population may have significant effects on other populations.
Therefore to apply the mathematical findings in real worlds, one must consider the interactions among the species in the ecosystem.
It seems that some results of this paper can be extended to multiple interacting species.
% under more general settings.
 For example,  a verification theorem like \thmref{thm-verification} can be established using almost the  same argument.
But one may no longer get  closed-form   value functions and   optimal controls by solving the corresponding quasi-variational inequalities.
 Also, we may be able to show that the value function is a viscosity solution of the quasi-variational inequalities but it
 %seems hard to
 is not immediate to
  identify condition(s) under which the value function is continuous.

In view of \cite[Section 4]{Hauss-S-I},
within the same framework and the same optimality criterion considered in this paper,
we may consider a more general problem  where the controlled state process is given by
\bed d\hat X(t) = b(\hat X(t),\al(t))dt + \sigma(\hat X(t),\al(t))dw(t)- \gamma(\hat X(t-),\al(t-))dZ(t),\eed
where $\hat X, Z\in \rr^n$,  $b,\sigma,\gamma$ are suitable functions with appropriate dimensions,
 and $w$ is an $n$-dimensional standard Brownian motion.
 Since every process of bounded variation  can be written as a difference of two nondecreasing processes,
   the control space can be enlarged by allowing
the singular control
  $Z$   to be an adapted  process with bounded variation.

A number of other questions deserve further investigations.
 In particular, in many practical situations, it is virtually impossible to obtain the explicit form of the value function and an optimal control by solving
\eqref{quasi-vari}.
Therefore a viable alternative is to employ numerical approximations. The controlled
Markov chain approximation method developed in \cite{Kushner-D} seems promising.
We may also consider relaxed control, under which we may achieve the optimal value
with an optimal control.
 Another problem of great interests is to consider the case when the random environment or the Markov chain $\al$ is unobservable.

%\bibliography{refs}\end{document}

\begin{thebibliography}{ }
\parskip=0pt
\bibitem{Alvarez}
L.H.R. Alvarez.
\newblock Singular stochastic control in the presence of a state-dependent
  yield structure.
\newblock {\em Stochastic Process. Appl.}, 86:323--343, 2000.

\bibitem{A-Shepp}
L.H.R. Alvarez and L.A. Shepp.
\newblock Optimal harvesting of stochastically fluctuating populations.
\newblock {\em J. Math. Biol.}, 37:155--177, 1998.

\bibitem{Arnold-H-S-79}
L.~Arnold, W.~Horsthemke, and J.W. Stucki.
\newblock The influence of external real and white noise on the
  {L}otka-{V}olterra model.
\newblock {\em Biomedical J.}, 21:451--471, 1979.

\bibitem{Brauman-02}
C.A. Brauman.
\newblock Variable effort harvesting models in random environments:
  generalization to density-dependent noise intensities.
\newblock {\em Mathematical Biosciences}, 177 \& 178:229--245, 2002.

\bibitem{Choulli}
T.~Choulli, M.~Taksar, and X.Y. Zhou.
\newblock A diffusion model for optimal dividend distribution for a company
  with constraints on risk control.
\newblock {\em SIAM J. Control Optim.}, 41(6):1946–--1979, 2003.

\bibitem{Cohen-90}
J.E. Cohen, T.~{\L}uczak, C.M. Newman, and Z.-M. Zhou.
\newblock Stochastic structure and nonlinear dynamics of food webs: qualitative
  stability in a {L}otka-{V}olterra cascade model.
\newblock {\em Proc. R. Soc. Lond. B}, 240:607--627, 1990.

\bibitem{Du-Sam-06}
N.H. Du and V.H. Sam.
\newblock Dynamics of a stochastic {L}otka-{V}olterra model perturbed by white
  noise.
\newblock {\em J. Math. Anal. Appl.}, 324:82--97, 2006.

\bibitem{FlemingS}
W.H. Fleming and H.M. Soner.
\newblock {\em Controlled Markov Processes and Viscosity Solutions}, volume~25
  of {\em Stochastic Modelling and Applied Probability}.
\newblock Springer-Verlag, New York, NY, second edition, 2006.

\bibitem{Guo-Z}
X.~Guo and Q.~Zhang.
\newblock Closed-form solutions for perpetual american put options with regime
  switching.
\newblock {\em SIAM J. Appl. Math.}, 64:2034--2049, 2004.

\bibitem{Guo-Z05}
X.~Guo and Q.~Zhang.
\newblock Optimal selling rules in a regime switching model.
\newblock {\em IEEE Transactions on Automatic Control}, 50:1450--1455, 2005.

\bibitem{Hauss-S-I}
U.~Haussmann and W.~Suo.
\newblock Singular optimal stochastic controls {I}: {E}xistence.
\newblock {\em SIAM J. Control Optim.}, 33(3):916--936, 1995.

\bibitem{Hauss-S}
U.~Haussmann and W.~Suo.
\newblock Singular optimal stochastic controls {II}: Dynamic programming.
\newblock {\em SIAM J. Control Optim.}, 33(3):937--959, 1995.

\bibitem{Jeffries76}
C.~Jeffries.
\newblock Stability of predation ecosystem models.
\newblock {\em Ecology}, 57:1321--1325, 1976.

\bibitem{Karatzas-S}
I.~Karatzas and S.E. Shreve.
\newblock {\em Brownian Motion and Stochastic Calculus}.
\newblock Springer, New York, second edition, 1991.

\bibitem{K}
R.Z. Khasminskii.
\newblock {\em Stochastic Stability of Differential Equations}.
\newblock Sijthoff and Noordhoff, Alphen aan den Rijn, Netherlands, 1980.

\bibitem{Kushner-D}
H.J. Kushner and P.~Dupuis.
\newblock {\em Numerical Methods for Stochstic Control Problems in Continuous
  Time}, volume~24 of {\em Stochastic Modelling and Applied Probability}.
\newblock Springer, New York, second edition, 2001.

\bibitem{Lande}
R.~Lande, S.~Engen, and B.~S{\ae}ther.
\newblock Optimal harvesting of fluctuating populations with a risk of
  extinction.
\newblock {\em The American Naturalist}, 145:728--745, 1995.

\bibitem{L-Oksendal}
E.M. Lungu and B.~{\O}ksendal.
\newblock Optimal harvesting from a population in a stochastic crowded
  environment.
\newblock {\em Math. Biosci.}, 145:47--75, 1997.

\bibitem{Lungu-O}
E.M. Lungu and B.~{\O}ksendal.
\newblock Optimal harvesting from interacting populations in a stochastic
  environment.
\newblock {\em Bernoulli}, 7:527--539, 2001.

\bibitem{Luo-Mao07}
Q.~Luo and X.~Mao.
\newblock Stochastic population dynamics under regime switching.
\newblock {\em J. Math. Anal. Appl.}, 334:69--84, 2007.

\bibitem{MSXZ08a}
J.~Ma, Q.S. Song, J.~Xu, and J.~Zhang.
\newblock Impulse control and optimal portfolio selection with general
  transaction cost.
\newblock 2008.
\newblock preprint.

\bibitem{Medina-R}
C.E. Medina-Reyna.
\newblock Growth and emigration of white shrimp, {L}itopenaeus vannamei, in the
  {M}ar {M}uerto {L}agoon, {S}outhern {M}exico.
\newblock {\em Naga, The ICLARM Quarterly}, 24:30--34, 2001.

\bibitem{Miller-V}
R.A. Miller and K.~Voltaire.
\newblock A stochastic analysis of the three paradigm.
\newblock {\em J. Econ. Dyn. Control}, 6:371--386, 1983.

\bibitem{Oksendal}
B.~{\O}ksendal.
\newblock {\em Stochastic differential equations, An introduction with
  applications}.
\newblock Springer-Verlag, Berlin, 6th edition, 2003.

\bibitem{OS02}
B.~{\O}ksendal and Agn{\`e}s Sulem.
\newblock Optimal consumption and portfolio with both fixed and proportional
  transaction costs.
\newblock {\em SIAM J. Control Optim.}, 40:1765--1790, 2002.

\bibitem{Otuma-O}
M.O. Otuma and I.I. Osakwe.
\newblock Assessment of the reproductive performance and post weaning growth of
  crossbred goat in derived {G}uinea {S}avanna {Z}one.
\newblock {\em Research Journal of Animal Sciences}, 2:87--91, 2008.

\bibitem{Paulsen-03}
J.~Paulsen.
\newblock Optimal dividend payouts for diffusions with solvency constraints.
\newblock {\em Finance Stoch.}, 7:457--–473, 2003.

\bibitem{Pha09}
H.~Pham.
\newblock {\em Continuous-time stochastic control and optimization with
  financial applications}, volume~61 of {\em Stochastic Modelling and Applied
  Probability}.
\newblock Springer-Verlag, Berlin, 2009.

\bibitem{Ryan-Hanson}
D.~Ryan and F.B. Hanson.
\newblock Optimal harvesting of a logistic population in an environment with
  stochastic jumps.
\newblock {\em J. Math. Biol.}, 24:259--277, 1986.

\bibitem{Sayre}
N.F. Sayre.
\newblock The genesis, history, and limits of carrying capacity.
\newblock {\em Annals of the Association of American Geographers}, 98:120--134,
  2008.

\bibitem{Slatkin-78}
M.~Slatkin.
\newblock The dynamics of a population in a {M}arkovian environment.
\newblock {\em Ecology}, 59:249--256, 1978.

\bibitem{Taksar-Z}
M.~Taksar and X.~Zeng.
\newblock On maximizing {CRRA} utility in regime switching markets with random
  endowment.
\newblock {\em SIAM J. Control Optim.}, 48:2984--3002, 2009/10.

\bibitem{Weerasinghe}
A.~Weerasinghe and A.~Mandelbaum.
\newblock {\em Abandonment vs. blocking in many-server queues: asymptotic
  optimility in the {QED} regime}, 2009.
\newblock preprint.

\bibitem{Wein90}
L.M. Wein.
\newblock Optimal control of a two station brownian network.
\newblock {\em Math. Oper. Res.}, 15:215--242, 1990.

\bibitem{YZ-10}
G.~Yin and C.~Zhu.
\newblock {\em Hybrid Switching Diffusions: Properties and Applications},
  volume~63 of {\em Stochastic Modelling and Applied Probability}.
\newblock Springer, New York, 2010.

\bibitem{YongZ}
J.~Yong and X.Y. Zhou.
\newblock {\em Stochastic Controls: Hamiltonian Systems and HJB Equations},
  volume~43 of {\em Stochastic Modelling and Applied Probability}.
\newblock Springer-Verlag, New York, 1999.

\bibitem{Zhou-Yin}
X.Y. Zhou and G.~Yin.
\newblock Markowitz's mean variance portfolio selection with regime switching:
  {A} continuous-time model.
\newblock {\em SIAM J. Control Optim.}, 42:1466--–1482, 2003.

\bibitem{ZY-09a}
C.~Zhu and G.~Yin.
\newblock On competitive {L}otka-{V}olterra model in random environments.
\newblock {\em J. Math. Anal. Appl.}, 357:154--170, 2009.

\end{thebibliography}

\def\cprime{$'$}

\end{document}